\documentclass[11pt,a4paper, fleqn]{article}
\usepackage[utf8]{inputenc}
\usepackage[T2A,T1]{fontenc}
\usepackage[english,russian]{babel}
\usepackage{epsfig,array}

\usepackage{mathtools}
\mathtoolsset{showonlyrefs}

\usepackage{indentfirst} % Красная строка
\usepackage[14pt]{extsizes}
\usepackage{amssymb}
\usepackage{amsmath}
\usepackage{amsthm}
\usepackage{pifont}
\usepackage{graphicx}
\usepackage{wrapfig}
\usepackage[colorinlistoftodos]{todonotes}
\usepackage[colorlinks=true, allcolors=blue]{hyperref}

\usepackage{algorithm}
\usepackage{algorithmic}
\usepackage{titlesec}
\titleformat{\section}
  {\large\bfseries}{\thesection.}{1em}{}
\titleformat{\subsection}
  {\centering\normalfont\normalsize\bfseries\itshape}{\thesubsection.}{1em}{}
\addto\captionsrussian{}
\allowdisplaybreaks

\newtheorem{theorem}{Теорема}
\newtheorem{corollary}{Следствие}[theorem]

\newtheorem{assumption}{Предположение}

\def\UDK#1{{\leftline{УДК {#1}}}}

\usepackage{geometry}
\def \PP {\mathbb P}
\def \EE {\mathbb E}
\def \R  {\mathbb R}

\def\E{{\mathbb{E}}}

\usepackage[labelsep=period]{caption}
\usepackage{subcaption}
\usepackage[backend=biber,bibencoding=utf8,sorting=none,style=gost-numeric,language=autobib,autolang=other,clearlang=true,sortcites=true,doi=false,isbn=false,date=year]{biblatex}
\newtheorem{lemma}{Лемма}
\newtheorem{definition}{Определение}

\begin{document}
\renewcommand{\abstractname}{\vspace{-\baselineskip}}

$$
\\\\
$$

{\it \UDK{519.85}}

\begin{center}
\textbf{Унифицированный анализ методов решения вариационных неравенств: редукция дисперсии, сэмплирование, квантизация и покомпонентный спуск}\footnote{
Исследование выполнено при поддержке Министерства науки и высшего образования Российской Федерации (госзадание) №075-00337-20-03, номер проекта 0714-2020-0005.
}

\textbf{\copyright~2021~г. 
А.\,Н. Безносиков$^{1}$,
А.\,В.~Гасников$^{1,2,3}$,
К.\,Э. Зайнуллина$^{1}$,
А.\,Ю. Масловский$^{1}$,
Д.\,А.~Пасечнюк$^{1,2}$,
}% В ЖВМ алфавитный порядок авторов. Так у них принято!

\textit{
$^{1}$141701 Московская обл., Долгопрудный, Институтский пер. 9, 
Московский физико-технический институт (национальный исследовательский университет), Россия \\
$^{2}$127051 Москва, Большой Каретный переулок, д.19 стр. 1, Институт проблем передачи информации им. А.А. Харкевича, Россия\\
$^{3}$385000 Республика Адыгея,  Майкоп, ул. Первомайская, д. 208 Кавказский математический центр Адыгейского государственного университета, Россия 
% $^{3}$Институт прикладного анализа и стохастики им. К. Вейерштрасса, Берлин, Германия \\
% $^{4}$109028 Москва, Покровский бульвар 11, 
% Национальный исследовательский университет <<Высшая школа экономики>>,
% Россия \\
% {$^*$ e-mail: asivanova@hse.ru}
}

{Поступила в редакцию 27.01.2022 г.\\
Переработанный вариант  00.00.2022 г.\\Принята к публикации 00.00.2022 г.}

\end{center}

\begin{abstract}
\noindent В данной статье предлагается унифицированный анализ методов для такого широкого класса задач, как вариационные неравенства, который в качестве частных случаев включает в себя задачи минимизации и задачи нахождения седловой точки. Предлагаемый анализ развивается на основе экстраградиентного метода, являющегося стандартным для решения вариационных неравенств. Рассматриваются монотонный и сильно монотонный случаи, которые соответствуют выпукло-вогнутым и сильно-выпукло-сильно-вогнутым задачам нахождения седловой точки. Теоретический анализ основан на параметризованных предположениях для итераций экстраградиентного метода. Следовательно, он может служить прочной основой для объединения уже существующих методов различных типов, а также для создания новых алгоритмов. В частности, чтобы показать это, мы разрабатываем некоторые новые надежные методы, в том числе метод с квантизацией, покомпонентный метод, распределенные рандомизированные локальные методы и другие. Большинство из упомянутых подходов прежде никогда не рассматривались в общности вариационных неравенств и применялись лишь для задач минимизации. Стабильность новых методов подтверждается предоставляемыми численными экспериментами по обучению моделей GAN. Библ.35. Фиг.3. Табл.1.

\textbf{Ключевые слова}: Экстраградиентный метод, стохастические вариационные неравенства, квантизация, редукция дисперсии.
\end{abstract}

%==========================================================================================
\section{Введение}

Основная постановка задачи, рассматриваемая в этой статье, является задачей решения вариационного неравенства (ВН), и имеет следующий вид:
\begin{equation} \label{VI}
    \text{найти} \quad z^* \in \mathcal{Z} \quad \text{такое, что} \quad \langle F(z^*), z - z^* \rangle + h(z) - h(z^*) \geq 0, \quad \forall z \in \mathcal{Z},
\end{equation}
где $h: \R^d \to \R \cup \{ + \infty\}$~--- подходящая полунепрерывная снизу выпуклая функция, $F: \R^d \to \R^d $~--- оператор, $\mathcal{Z}$~--- непустое замкнутое выпуклое подмножество $\R^d$. 

Приведем сразу некоторые классические примеры задач, которые могут быть представлены в виде \eqref{VI}:

\textbf{Задача минимизации.} Положим, что $F(x)$ является градиентом некоторой функции $f(x)$, а $h(x) = \delta_{\mathcal{X}}(x)$~--- индикаторная функция множества $\mathcal{X}$. В частности, при $\mathcal{X} = \R^d$, $x^* \in \mathcal{X}$ является решением \eqref{VI} тогда и только тогда, когда $\nabla f(x^*) = 0$. В случае выпуклой функции, ему соответствует глобальный экстремум.

\textbf{Задача нахождения седловой точки (СТ).}
Рассмотрим следующую минимаксную задачу:
\begin{equation}
    \label{minmax}
    \min_{x \in \mathcal{X}} \max_{y \in \mathcal{Y}} f(x, y).
\end{equation}
Эта задача может быть также представлена в виде \eqref{VI}. Достаточно выбрать $F$ и $g$ следующим образом:
\begin{equation*}
    F(z) = \left(
\begin{array}{c}
\nabla_x f(x,y) \\
-\nabla_y f(x,y)
\end{array}
\right), \quad \quad
h(z) = \delta_{\mathcal{X}}(x) + \delta_{\mathcal{Y}}(y).
\end{equation*}
Как и в случае задачи минимизации, неравенство \eqref{VI} является необходимым условием оптимальности. В частности, если $f(x,y)$ выпукло-вогнута, и $\mathcal{X} = \R^{d_x}$, $\mathcal{Y} = \R^{d_y}$, это условие является также и достаточным, причём решение $(x^*, y^*)$ вариационного неравенства тогда является глобальной седловой точкой:
\begin{equation*}
    f(x^*, y) \leq f(x, y) \leq f(x, y^*), \quad \forall x \in \mathcal{X}, y \in \mathcal{Y}.
\end{equation*}

Эти примеры показывают, что класс задач ВН достаточно широк. В частности, задачи минимизации могут рассматриваться в общности ВН. Но обычно предпочитают проводить анализ задачи минимизации отдельно и независимо. На данный момент, анализ задач минимизации разработан гораздо шире и глубже, чем для задач ВН. В первую очередь это связано со сложностью задачи ВН: многие техники из минимизации ещё не перенесены или не могут быть перенесены на ВН. Поэтому прежде всего решение ВН интересно с точки зрения нахождения СТ.

Задачи СТ имеют множество практических приложений. Как, например, хорошо известные классические примеры из теории игр и оптимального управления \cite{facchinei2007finite}. В последние годы задачи СТ стали популярными и в нескольких других отношениях. Можно отметить ветвь работ, посвященных решению негладких задач путем их переформулирования в виде задачи СТ \cite{nesterov2005smooth, nemirovski2004prox}, а также применение таких подходов к задачам обработки изображений \cite{chambolle2011first, esser2010general}. Однако в первую очередь интересны приложения седловых задач в машинном обучении. Конечно, прежде всего здесь стоит упомянуть о GAN. В классической формулировке \cite{goodfellow2014generative} обучение этих моделей является минимаксной задачей.

Довольно большая часть прикладных задач, в том числе задачи машинного обучения, являются стохастическими, поэтому естественно сосредоточиться на случае, когда невыгодно (или даже невозможно) вычислить полное значение градиента и когда вместо этого используются некоторые стохастические оценки. Например, для функции из \eqref{minmax}, имеющей следующий вид: $f(x, y) =  \EE_{p_{x} \sim D_x, p_{y} \sim D_y} \left[f_{p_{x}, p_{y}}(x, y)\right]$. В частности, в случае GAN, $f$ есть функция потерь, а переменные $x$ и $y$ интерпретируются как относящиеся к двум моделям: $x$~--- параметры дискриминатора, а $y$~--- генератора, $p_{x}$~--- обучающий пример из реального набора данных, $p_{y}$~--- случайный вектор, который генератор использует для создания поддельных копий реального набора данных. Стандартное предположение статистического обучения состоит в том, что распределение данных $D_x$ неизвестно, а потому полный градиент $\nabla_x f(x, y)$ не может быть вычислен, тогда как можно легко вычислить градиент для некоторых отдельных данных. Возвращаясь к основной проблеме \eqref{VI} этой статьи, мы интерпретируем сказанное выше следующим образом: предположим, что мы не имеем доступа к ``честному'' $F(z)$, а только к некоторой несмещенной стохастической оценке $F(z, \xi)$:
\begin{equation}
    \label{stoc}
    F(z) = \EE_{\xi \sim P} \left[F (z,\xi)\right].
\end{equation}

В данной работе нас также будет интересовать другая стохастическая постановка задачи \eqref{VI}, это тот случай, когда значение $F$~--- это среднее значение большого числа операторов:
\begin{equation}
    \label{MK}
    F(z) = \frac{1}{M}\sum\limits_{m=1}^M F_m(z).
\end{equation}
Такие постановки возникают в результате применения подхода интегрирования  Монте--Карло. Например, пусть для \eqref{minmax} имеется $M$ частей набора данных, тогда можно вычислить градиенты $\nabla_x f_m(x, y)$, $\nabla_y f_m(x, y)$ для каждой из этих частей, тогда как вычисление полных градиентов $\nabla_x f(x, y)$, $\nabla_y f(x, y)$ будет очень дорогим и потребует много времени. Отсюда кажется естественным выбирать случайный индекс $m$ части набора данных на каждой итерации и учитывать градиенты только на ней. Этот подход обычно используется при практическом решении задач обучения.

Приведенный выше взгляд на \eqref{MK} справедлив в случае, когда у нас есть только одно устройство для вычислений. Однако в случае распределенного обучения можно просто обмениваться данными между устройствами, и тогда каждому устройству будет соответствовать свой $F_m$. В такой постановке задачу также можно рассматривать в виде \eqref{MK}, но мы предпочтём переписать ее следующим образом:
\begin{equation}
    \label{distr}
    F(Z) = \Phi(Z) + \lambda \cdot (Z - \bar Z ),
\end{equation}
где вектора $\Phi(Z) = [F^T_1(z_1), \ldots, F^T_M(z_M)] \in \R^{Md}$,  $Z = [z_1^T, \ldots, z_M^T] \in \R^{Md}$ и $\bar Z = [\bar z^T, \ldots, \bar z^T] \in \R^{Md}$ для $ \bar z = \frac{1}{M} \sum_{m=1}^M  z_m$. Легко проверить, что задача минимизации $\min_{x_1, \ldots x_M} \Big[\sum_{m=1}^M f_m(x_m) + \frac{\lambda}{2} \sum_{m=1}^M \| x_m - \bar x \|^2\Big]$ соответствует ВН с оператором \eqref{distr}. Аналогично, задача СТ $$\min_{x_1, \ldots x_M} \max_{y_1, \ldots y_M} \Big[\sum_{m=1}^M f_m(x_m, y_m) + \frac{\lambda}{2} \sum_{m=1}^M \| x_m - \bar x \|^2 - \frac{\lambda}{2} \sum_{m=1}^M \| y_m - \bar y \|^2\Big]$$ также соответствует ВН с оператором \eqref{distr} (как построить ВН из задачи минимизации и задачи нахождения седловой точки описано выше). Разобравшись в интуиции, переходим к интерпретации: $z_1, \ldots z_M$~--- локальные модели на каждом устройстве, $\bar z$~--- глобальная модель, полученная усреднением всех локальных моделей. Понятно, что нужно выбрать параметр регуляризации $\lambda$ достаточно большим, чтобы решения задач \eqref{MK} и \eqref{distr} были примерно одинаковыми: $z_m \approx \bar z$. Но имеет ли смысл брать $\lambda $ малым (например, когда $\lambda = 0$, мы имеем просто локальные модели)? Оказывается, да: такая постановка задачи с переменной $\lambda$, значение которой может варьироваться, породила целую тенденцию в персонализированном федеративном обучении~--- федеративное обучение со смешиванием \cite{hanzely2020federated,hanzely2020lower}. 
% Problem \eqref{distr} can be generalized to the case of decentralized learning by replacing the regularizer: $\lambda(Z - \bar Z ) \to \lambda WZ$. Such a problem also arises in federated learning, namely in multi-task learning \cite{wang2018distributed,smith2017federated}. We describe it in more detail in the appropriate section of the Appendix. Setting \eqref{distr} and its decentralized generalization is the last formulation of \eqref{VI} to be explored by our unified analysis.
Постановка \eqref{distr} будет последней из формулировок задачи \eqref{VI}, которые будут исследоваться в рамках предлагаемого унифицированного анализа.

\subsection{Экстраградиентный метод}

Стохастический градиентный спуск по-прежнему является основным методом минимизации, он используется для большого количества задач машинного обучения, несмотря на наличие более современных методов. Основным же методом решения гладких вариационных неравенств является экстраградиентный метод \cite{Korpelevich1976TheEM}. В простом виде методы такого рода записываются следующим образом:
\begin{equation}
    \label{eg}
    z^{k+1/2} =\text{prox}_{\gamma h} ( z^k - \gamma g^{k}), \quad \quad z^{k+1} = \text{prox}_{\gamma h} (z^k - \gamma g^{k+1/2}),
\end{equation}
где $\text{prox}_{\gamma h} (z) = \arg\min_x \left\{ \gamma h(x) + \frac{1}{2} \| z - x\|^2\right\}$. 
В классическом, детерминированном, варианте $g^k = F(z^k)$ и $g^{k+1/2} = F(z^{k+1/2})$. Этот метод оптимален с точностью до численных констант для гладких монотонных и сильно монотонных вариационных неравенств \cite{juditsky2011solving}. Причем на практике этот метод показывает себя лучше, чем спуск с итерацией обычного вида $z^{k+1} = z^k - \gamma g^{k}$. Более того, известно, что метод без дополнительного шага расходится для наиболее распространенных билинейных задач. Следовательно, вычисление $z^{k + 1/2}$ является ключевым.
В этой работе, для нашего унифицированного анализа, мы используем метод \eqref{eg} с несколько более сложной структурой:
\begin{align}
    \label{eg_mod}
    \bar z^k &= \tau z^k + (1 - \tau) w^k, \nonumber\\
    z^{k+1/2} &= \text{prox}_{\gamma h} (\bar z^k - \gamma g^{k}),\\
     z^{k+1} &= \text{prox}_{\gamma h}(\bar  z^k - \gamma g^{k+1/2}), \nonumber
\end{align}
где $w^{k+1} = z^{k+1}$ с вероятностью $(1 - \tau)$ или $w^k$, иначе. Легко увидеть, что при $\tau = 0$ верно $w^k = z^k$, и метод \eqref{eg_mod} в точности соответствует \eqref{eg}. Метод \eqref{eg_mod} был впервые предложен в \cite{alacaoglu2021stochastic}. 

\subsection{От минимизации к ВН}

SGD используется для решения задач минимизации с середины прошлого века \cite{10.1214/aoms/1177729586} и с того времени был расширен огромным количеством различных модификаций. Это методы уменьшения дисперсии \cite{johnson2013accelerating}, квантования \cite{alistarh2017qsgd} и координатные методы \cite{hanzely2018sega} и т.д. Различные варианты SGD см. в \cite{gorbunov2020unified}. В отличие от задач минимизации, вариационные неравенства и седловые задачи не имеют такого широкого набора теоретических результатов, хотя базовый метод для задачи ВН и более сложен и дает широкий простор для творчества. Но в то же время развитие тех же идей, что и для задач минимизации, для ВН происходит значительно медленнее, в том числе и из-за того, что ВН более общи и сложны в теоретическом анализе. Далее мы перечислим основные достижения в области решения ВН касательно конструирования методов подобных уже существующим методам минимизации.

$\bullet$ \textbf{Базовые методы.} Как отмечалось ранее, базовым методом решения ВН является экстраградиентный метод \cite{Korpelevich1976TheEM}; еще более общая его версия называется Mirror Prox \cite{nemirovski2004prox}. Анализ в стохастическом случае с ограниченной дисперсией шума описан в \cite{juditsky2011solving}. 
Стоит обратить внимание на интересные модификации этого метода: с одним дополнительным вызовом оракула \cite{hsieh2019convergence} и с повторный вызовом \cite{mishchenko2020revisiting}. Можно также причислить следующие классические методы, отличающиеся по структуре от экстраградиентного: \cite{doi:10.1137/S0363012998338806,nesterov2007dual}.

$\bullet$ \textbf{Редукция дисперсии.} Направление разработки методов редукции дисперсии для задач ВН и СТ развивалось начиная с работы \cite{NIPS2016_1aa48fc4}, где представлен метод, основанный на сильно выпукло-сильно вогнутых седлах (сильно-монотонных ВН). Также в \cite{chavdarova2019reducing} был предложен метод для сильно монотонных ВН. Наконец, стоит выделить работу \cite{alacaoglu2021stochastic}, которая пересекается с прошлыми результатами или повторяет их, предоставляя методы редукции дисперсии для монотонных и сильно монотонных ВН. Отметим, что приведенные выше методы сильно отличаются друг от друга, более того, они далеки от классической редукции дисперсии для задач минимизации.

$\bullet$ \textbf{Покомпонентные и квантизованные методы.} Покомпонентные методы для задач СТ и ВН изучены не слишком хорошо. Можно выделить работы, посвященные конкретным покомпонентным методам для каких-то определенных классов задач СТ \cite{sidford2018coordinate, carmon2020coordinate}, а также работы по безградиентным методам \cite{sadiev2020zeroth}. Методы же с квантизацией, специализированные для задач СТ, до сих пор совсем не предлагались.

$\bullet$ \textbf{Локальные методы.} Направление разработки локальных методов изучено совсем не в той мере \cite{deng2021local,beznosikov2021distributed}, как это сделано для задач минимизации. В данной работе мы обращаем внимание не на детерминированные методы типа Local SGD, а на рандомизированные, которые могут быть применены для решения задачи \eqref{distr}, как например описанные в работе \cite{hanzely2020federated}. Суть этих методов в том, что мы с некоторой вероятностью вызываем и делаем шаг только по оракулу для $\Phi(Z)$, иначе обращаемся к $\lambda(Z - \bar Z)$. Вызов $\Phi (Z)$ соответствует локальной итерации, а вызов $\lambda (Z - \bar Z)$ -- коммуникации. 

\subsection{Наш вклад}

$\bullet$ \textbf{Унифицированный анализ.} В данной работе предлагается унифицированный теоретический анализ для методов типа \eqref{eg} и \eqref{eg_mod} в сильно монотонном и монотонном случаях. Анализ основан на параметризованных предположениях, поэтому он позволяет легко конструировать и анализировать огромное количество новых методов.

$\bullet$ \textbf{Улучшенные оценки для существующих методов.} В исходном анализе {\tt Past ES} из работы \cite{hsieh2019convergence} в стохастическом сильно монотонном случае достигается сублинейная скорость сходимости, тогда как мы можем гарантировать линейную скорость в детерминированном члене и сублинейную в стохастическом члене. Более того, мы предоставляем оценки для {\tt Past ES} и в стохастическом монотонном случае. Кроме того, мы покрываем результаты для некоторых других существующих методов или немного их обобщаем.

$\bullet$ \textbf{5 новых методов.} Но более важным, чем обобщение других результатов, является получение новых надежных методов. В отличие от работы \cite{gorbunov2020unified}, к моменту написания которой большинство анализируемых там методов уже были описаны в других работах, в нашей работе более половины методов являются новыми. Это покомпонентные методы, методы с квантизацией, методы с сэмплированием по важности, локальные рандомизированные методы.

$\bullet$ \textbf{Эксперименты.} Мы предоставляем сравнение методов на примере практической задачи, в котором демонстрируем, что предлагаемые нами новые методы могут превосходить по эффективности существующие. Эксперименты проводятся на искусственной билинейной задаче и, в некоторых случаях, на GAN.

%==============================================================================================
\section{Полученные результаты}\label{main_results}
Для начала введем основные определения. Мы используем аннотацию $\langle x,y \rangle := \sum_{i=1}^n x_i y_i$ для скалярного произведения векторов $x,y\in\R^d$. Оно порождает $\ell_2$-норму в пространстве $\R^d$ в таком виде: $\|x\| := \sqrt{\langle x, x \rangle}$.

\subsection{Основные предположения}

Далее нам потребуются два основных предположения. Наше первое предположение относится к монотонности оператора $F$ из \eqref{VI} и сильной выпуклости $g$. В частности, мы рассматриваем строго монотонный и монотонный случаи. С точки зрения задачи поиска седловой точки это соответствует сильно выпукло-сильно вогнутому и выпукло-вогнутому случаям.

\begin{assumption}\label{as1}
\textbf{(СМ) Сильная монотонность/сильная выпуклость.} Существуют неотрицательные $\mu_F, \mu_h$ такие, что $\mu_h + \mu_F > 0$  и следующие неравенства верны для всех $z_1, z_2 \in \mathcal{Z}$:
\begin{eqnarray}
    \label{sm1}
    \langle F(z_1) - F(z_2), z_1 - z_2 \rangle &\geq& \mu_F\|z_1-z_2\|^2 \\
    \label{sm2}
    h(z_1) - h(z_2) - \langle \nabla h(z_2), z_1 - z_2\rangle &\geq&  \frac{\mu_h}{2} \|z_1 - z_2 \|^2.
\end{eqnarray}

\textbf{(M) Монотонность/выпуклость.} Для всех $z_1, z_2 \in \mathcal{Z}$:
\begin{align}
    \label{eq_m}
    \langle F(z_1) - F(z_2), z_1 - z_2 \rangle \geq 0, \quad
    h(z_1) - h(z_2) - \langle \nabla h(z_2), z_1 - z_2\rangle \geq  0.
\end{align}

\end{assumption}

% The first two cases of Assumption \eqref{as1} are more than standard. The third (non-monotone case) considers a certain "good" class non-monotone operators (non-convex-non-concave saddles). For more information see,  for example, \cite{dang2015convergence}. (?)

Первые два случая в Предположении \eqref{as1} более чем стандартны. Третий же (не монотонны случай) рассматривает некоторый класс ``хороших'' немонотонных операторов (не выпуклых-не вогнутых сёдел). Подробности см., например, в \cite{dang2015convergence}.

Второе предположение ключевое и позволяет нам рассматривать разные методы для решения ВН в унифицированном виде. Суть этого предположения проста, аналогично с \cite{gorbunov2020unified}, мы вводим неравенства для основных членов, которые необходимо оценить при анализе экстра-градиентных методов.
 
\begin{assumption} \label{as2}
Пусть последовательности $\{z^k\}$ и $\{w^k\}$ были получены в результате случайных итераций \eqref{eg_mod}. Предположим, что стохастические операторы $g^{k+1/2}$ несмещены для всех $k$
\begin{align}
    \label{a0}
    \EE\left[g^{k+1/2} ~| ~z^{k+1/2} \right] = F(z^{k+1/2}),
\end{align}
Далее, предположим, что существуют неотрицательные $A,B,C,E,D_1$,  $D_2, D_3$ и $\rho \in (0;1]$ и рандомная последовательность $\{\sigma^k\}$ (может быть нулевой), что выполняются следующие неравенства:
%\textbf{Universal (for all cases):}
\begin{equation}
    \label{a1}
    \EE\left[\|g^{k+1/2} - g^k \|^2 \right] \leq A \EE\left[\| z^{k+1/2} -w^k \|^2\right] + B \EE\left[\sigma^2_k\right] + D_1,
\end{equation}
\begin{equation}
    \label{a2}
    \EE\left[\sigma^2_{k+1}\right] \leq (1 - \rho)\EE\left[\sigma^2_k\right]  + C \EE\left[\| z^{k+1/2} -w^k \|^2\right] + D_2,
\end{equation}
\begin{equation}
    \label{a3}
    \EE\left[\|g^{k+1/2} - F(z^{k+1/2}) \|^2 \right] \leq E \EE\left[\| z^{k+1/2} -w^k \|^2\right]  + D_3.
\end{equation}
\end{assumption}

Похожие неравенства могут быть найдены в анализе методов \eqref{eg} и \eqref{eg_mod}, см. \cite{juditsky2011solving,alacaoglu2021stochastic}.
Это первая работа, которая рассматривает вариационные неравенства в такой общности.

\subsection{Унифицированная теорема}
Мы готовы представить основной теоретический результат данной статьи. Для начала, введем функцию Ляпунова, с помощью которой будем анализировать сходимость:

\begin{equation}
    \label{lyap}
    V^k = \tau\| z^{k+1} - z^*\|^2 +  \| w^{k+1} - z^*\|^2 + T \gamma^2 \sigma^2_k,
\end{equation}
где константа $T > 0$. Этот критерий используется в сильно-монотонном случае.

Для монотонного случая используется другой критерий - функция зазора
\cite{nemirovski2004prox,juditsky2011solving}:
\begin{equation}
    \label{gap}
    \text{Gap} (z) = \max_{u \in \mathcal{C}} \left[ \langle F(u),  z - u  \rangle + h(z) - h(u) \right].
\end{equation}
Здесь максимум берется не по всему множеству $\mathcal{Z}$ (как в классической версии), а по $\mathcal{C}$ -- компактному подмножеству множества $\mathcal{Z}$.
Таким образом мы можем рассматривать неограниченные множества $\mathcal{Z}$. Это допустимо, так как такой вариант критерия верен, если решение $z^{*}$ лежит в  $\mathcal{C}$, см.  \cite{nesterov2007dual}.

\begin{theorem}
\label{unified}
Пусть выполнено Предположение 2. Тогда, если дополнительно выполняется одно из условий Предположения 1, верны следующие оценки
\begin{itemize}
    \item {для сильно-монотонного/сильно-выпуклого} случаев с \\ $\gamma \leq \min\left\{ \frac{\sqrt{1 - \tau}}{2\sqrt{2A + TC}}; \frac{1 - \tau}{4\left( \mu_F + \mu_h\right)}\right\}$ и $T \geq \frac{4B}{\rho}$:
    \begin{align*}
    \EE\left[V_{K}\right] 
    \leq &\max\left\{\left( 1 -\gamma \cdot \frac{\mu_F + \mu_h}{16}  \right)^{K-1} ; \left(1 - \frac{\rho}{2} \right)^{K-1} \right\} V_{0} \\
    &+ \frac{\gamma^2 (2D_1 + T D_2)}{\min\left\{\gamma \cdot \frac{\mu_F + \mu_h}{16}; \frac{\rho}{2}\right\}};
    \end{align*}
    \item {для монотонного/выпуклого} случаев с $\gamma \leq \frac{\sqrt{1 -\tau}}{2\sqrt{2A + TC + E}}$ и $T \geq \frac{2B}{\rho}$
    \begin{align*}
    \EE\left[ \text{Gap} (\bar z^K)\right] 
    \leq&  \frac{8\max_{u \in \mathcal{C}}\left[\| z^0 - u \|^2\right]  + 4\gamma^2 T \sigma^2_{0}}{\gamma K}  \\
    &+ \gamma  (7 D_1 + 3TD_2 + D_3),
    \end{align*}
    где $\bar z^{K} = \frac{1}{K}\sum\limits_{k=0}^{K-1} z^{k+1/2} $.
\end{itemize}
\end{theorem}

Метод \eqref{eg_mod} в условиях Предположения \ref{as2} сходится линейно в сильно-монотонном случае и сублинейно $(\sim \frac{1}{K})$ в монотонном случае до определенного радиуса осцилляции сходимости, который зависит от второго члена в Теореме \ref{unified}. Формально, можно добиться сходимости по этому члену за фиксированное число шагов. Для этого необходимо правильно выбрать шаг $\gamma$ \cite{stich2019unified}. Для некоторых методов мы используем это, чтобы получить классические оценки сходимости -- можно посмотреть об этом подробнее в следующем разделе. Оптимальный выбор шага для каждого метода находится в приложении, соответствующем этому методу.

Далее приведем доказательство теоремы  \ref{unified}. Для начала докажем лемму.

\begin{lemma} \label{lemm1}
 Пусть $h$ $\mu_h$ -- сильно выпуклая и $z^+ = \text{prox}_{\gamma h} (z)$. Тогда для всех $x \in \R^d$ справедливо следующее неравенство:
\begin{equation}
    \langle z^+ - z, x - z^+ \rangle \geq \gamma \left( h(z^+) - h(x) + \frac{\mu_h}{2} \|z^+ - x \|^2 \right).
\end{equation}
\end{lemma}
\textbf{Доказательство:} Мы используем $\gamma \mu$-сильно выпуклость функции $\gamma h$ \eqref{sm2}: 
\begin{align*}
    \gamma \left(h(x) - h(z^+)\right) - \langle \gamma \nabla h(z^+), x - z^+\rangle \geq  \frac{\gamma \mu_h}{2} \|x - z^+ \|^2.
\end{align*}
Вместе с определением $\text{prox}$ и необходимым условием оптимума: $\gamma \nabla h(z^+) = z - z^+$, это завершает доказательство.

Доказательство Теоремы \ref{unified}: \\
\textbf{Доказательство:} По Лемме \ref{lemm1} для $z^{k+1/2} = \text{prox}_{\gamma h} (\bar z^k - \gamma g^k)$ and $z^{k+1} = \text{prox}_{\gamma h} (\bar z^k - \gamma g^{k+1/2})$ для $x = u$ получаем
\begin{align*}
&\langle z^{k+1} - \bar z^k + \gamma g^{k+1/2} , u - z^{k+1}\rangle \geq \gamma \left(h(z^{k+1}) - h(u) + \frac{\mu_h}{2} \|z^{k+1} - u\|^2 \right),\\
&\langle z^{k+1/2} - \bar z^k + \gamma g^{k} , z^{k+1} - z^{k+1/2}\rangle \\
&\hspace{4.5cm} \geq \gamma \left(h(z^{k+1/2}) - h(z^{k+1}) + \frac{\mu_h}{2} \|z^{k+1} - z^{k+1/2}\|^2 \right).
\end{align*}
Далее суммируем два неравенства и производим некоторые перестановки:
\begin{align*}
&\langle z^{k+1} - \bar z^k, u - z^{k+1}\rangle + \langle z^{k+1/2} - \bar z^k, z^{k+1} - z^{k+1/2}\rangle \\ 
&+ \gamma \langle g^{k+1/2} - g^{k} , z^{k+1/2} - z^{k+1} \rangle + \gamma \langle g^{k+1/2}, u - z^{k+1/2} \rangle \\
&\hspace{1.5cm}\geq \gamma \left( h(z^{k+1/2}) - h(u)  + \frac{\mu_h}{2} \|z^{k+1} - z^{k+1/2}\|^2 + \frac{\mu_h}{2} \|z^{k+1} - u\|^2 \right).
\end{align*}
Умножая на 2 и используя определение $\bar z^k$ из \eqref{eg_mod}, получаем
\begin{align*}
&2 \tau \langle z^{k+1} - z^k, u - z^{k+1}\rangle + 2 (1 - \tau) \langle z^{k+1} - w^k, u - z^{k+1}\rangle \\
&+ 2 \tau \langle z^{k+1/2} - z^k, z^{k+1} - z^{k+1/2}\rangle + 2 (1 - \tau) \langle z^{k+1/2} - w^k, z^{k+1} - z^{k+1/2}\rangle \\
&+ 2\gamma \langle g^{k+1/2} - g^{k} , z^{k+1/2} - z^{k+1} \rangle + 2\gamma \langle g^{k+1/2}, u - z^{k+1/2} \rangle \\
&\hspace{1.4cm}\geq 2\gamma \left( h(z^{k+1/2}) - h(u)  + \frac{\mu_h}{2} \|z^{k+1} - z^{k+1/2}\|^2 + \frac{\mu_h}{2} \|z^{k+1} - u\|^2 \right).
\end{align*}
Для первой и второй строки используем выражение $2 \langle a,b \rangle = \| a + b \|^2 - \| a\|^2 - \| b \|^2$, и получаем
\begin{align*}
&\tau \left(\| z^k - u \|^2 -
\| z^{k+1} - z^k\|^2 -
\| z^{k+1} - u \|^2 \right) \\
&+
(1 - \tau) \left(\| w^k - u \|^2 -
\| z^{k+1} - w^k\|^2 -
\| z^{k+1} - u\|^2 \right) \\
&+
\tau (\|z^{k+1} - z^k\|^2 -
\| z^{k+1/2} - z^k \|^2 -
\|z^{k+1} - z^{k+1/2}\|^2)  \\
&+
(1 -\tau) (\|z^{k+1} - w^k\|^2 -
\| z^{k+1/2} - w^k \|^2 -
\|z^{k+1} - z^{k+1/2}\|^2)  \\
&+ 2\gamma \langle g^{k+1/2} - g^{k} , z^{k+1/2} - z^{k+1} \rangle + 2\gamma \langle g^{k+1/2}, u - z^{k+1/2} \rangle \\
&\hspace{1.4cm} \geq 2\gamma \left( h(z^{k+1/2}) - h(u)  + \frac{\mu_h}{2} \|z^{k+1} - z^{k+1/2}\|^2 + \frac{\mu_h}{2} \|z^{k+1} - u\|^2 \right).
\end{align*}
Небольшая перестановка дает
\begin{align*}
(1 &+\gamma \mu_h )\| z^{k+1} - u\|^2 \leq \tau \| z^k - u \|^2 + (1 - \tau) \| w^k - u \|^2 \\
&- \tau \| z^{k+1/2} - z^k\|^2 - (1 -\tau) \| z^{k+1/2} - w^k \|^2 \\
&- (1 + \gamma \mu_h ) \|z^{k+1} - z^{k+1/2}\|^2 \\
&+ 2\gamma \langle g^{k+1/2} - g^{k} , z^{k+1/2} - z^{k+1} \rangle - 2\gamma \langle g^{k+1/2},  z^{k+1/2} - u  \rangle  \\
&- 2\gamma \left( h(z^{k+1/2}) - h(u)\right).
\end{align*}
Из простого факта: $2 \langle a,b \rangle \leq  \eta \| a\|^2 +  \frac{1}{\eta}\| b \|^2$ с $a = g^{k+1/2} - g^{k}, b =z^{k+1/2} - z^{k+1}, \eta = 2\gamma$, следует
\begin{align}
\label{temp0}
(1 &+ \gamma \mu_h )\| z^{k+1} - u\|^2 \leq \tau \| z^k - u \|^2 + (1 -\tau) \| w^k - u \|^2 \notag\\
&- \tau \| z^{k+1/2} - z^k\|^2 - (1 -\tau) \| z^{k+1/2} - w^k \|^2 \notag\\ 
& - \left(\frac{1}{2} + \gamma \mu_h\right)  \|z^{k+1} - z^{k+1/2}\|^2 + 2\gamma^2 \| g^{k+1/2} - g^{k}\|^2 \notag\\ 
&- 2\gamma \langle g^{k+1/2},  z^{k+1/2} - u  \rangle - 2\gamma \left( h(z^{k+1/2}) - h(u)\right) .
\end{align}
Далее, рассмотрим разные случаи теоремы. Начнем с \textbf{сильно монотонного/выпуклого случая}. Заменим $u = z^*$, возьмем полное математическое ожидание и получим
\begin{align*}
(1 &+ \gamma \mu_h ) \EE\left[\| z^{k+1} - z^*\|^2 \right] \leq \tau \EE\left[\| z^k - z^* \|^2\right] + (1 -\tau) \EE\left[\| w^k - z^* \|^2\right] \\
&\hspace{0.4cm}- \tau \EE\left[\| z^{k+1/2} - z^k\|^2\right] - (1 -\tau) \EE\left[\| z^{k+1/2} - w^k \|^2\right] \\
&\hspace{0.4cm} - \left(\frac{1}{2} + \gamma \mu_h\right) \EE\left[\|z^{k+1} - z^{k+1/2}\|^2\right] + 2\gamma^2 \EE\left[\| g^{k+1/2} - g^{k}\|^2\right] \\
&\hspace{0.4cm} - 2\gamma \EE\left[\langle g^{k+1/2},  z^{k+1/2} - z^*  \rangle + h(z^{k+1/2}) - h(z^*)\right] \\
&= \tau \EE\left[\| z^k - z^* \|^2\right] + (1 -\tau) \EE\left[\| w^k - z^* \|^2\right]  \\
&\hspace{0.4cm} - (1 -\tau) \EE\left[\| z^{k+1/2} - w^k \|^2\right] \\
&\hspace{0.4cm} - \left(\frac{1}{2} + \gamma \mu_h\right)  \EE\left[\|z^{k+1} - z^{k+1/2}\|^2\right] + 2\gamma^2 \EE\left[\| g^{k+1/2} - g^{k}\|^2\right] \\
&\hspace{0.4cm} - 2\gamma \EE\left[\langle \EE\left[g^{k+1/2} ~|~ z^{k+1/2}\right],  z^{k+1/2} - z^*  \rangle + h(z^{k+1/2}) - h(z^*)\right].
\end{align*}
Далее применим Предположение \ref{as2}, а именно, \eqref{a0} и \eqref{a1}:
\begin{align*}
(1 &+\gamma  \mu_h ) \EE\left[\| z^{k+1} - z^*\|^2 \right] \leq \tau \EE\left[\| z^k - z^* \|^2\right] + (1 -\tau) \EE\left[\| w^k - z^* \|^2\right] \\
&\hspace{0.4cm}- \tau \EE\left[\| z^{k+1/2} - z^k\|^2\right] - (1 -\tau) \EE\left[\| z^{k+1/2} - w^k \|^2\right] \\
&\hspace{0.4cm} - \left(\frac{1}{2} + \gamma \mu_h\right) \EE\left[\|z^{k+1} - z^{k+1/2}\|^2\right] \\
&\hspace{0.4cm} + 2\gamma^2 \left( A \EE\left[\| z^{k+1/2} -w^k \|^2\right] + B \EE\left[\sigma^2_k\right] + D_1 \right)\\ 
&\hspace{0.4cm}- 2\gamma \EE\left[\langle F(z^{k+1/2}),  z^{k+1/2} - z^*  \rangle + h(z^{k+1/2}) - h(z^*)\right]  \\
&= \tau \EE\left[\| z^k - z^* \|^2\right] + \left(1 -\tau \right) \EE\left[\| w^k - z^* \|^2\right] - \tau \EE\left[\| z^{k+1/2} - z^k\|^2\right] \\
&\hspace{0.4cm} - \left(\frac{1}{2} + \gamma \mu_h\right) \EE\left[\|z^{k+1} - z^{k+1/2}\|^2\right] \\
&\hspace{0.4cm} - \left((1 -\tau) - 2\gamma^2 A \right) \EE\left[\| z^{k+1/2} - w^k \|^2\right] + 2\gamma^2 B \EE\left[\sigma^2_k\right] \\
&\hspace{0.4cm} - 2\gamma \EE\left[\langle F(z^{k+1/2}),  z^{k+1/2} - z^*  \rangle + h(z^{k+1/2}) - h(z^*)\right] + 2\gamma^2 D_1 .
\end{align*}
Свойство решения \eqref{VI} дает
\begin{align*}
(1 &+\gamma  \mu_h ) \EE\left[\| z^{k+1} - z^*\|^2 \right] \leq \tau \EE\left[\| z^k - z^* \|^2\right] + \left(1 -\tau \right) \EE\left[\| w^k - z^* \|^2\right] \\
&\hspace{0.4cm} - \tau \EE\left[\| z^{k+1/2} - z^k\|^2\right] - \left(\frac{1}{2} + \gamma \mu_h\right) \EE\left[\|z^{k+1} - z^{k+1/2}\|^2\right]  \\
&\hspace{0.4cm} - \left((1 -\tau) - 2\gamma^2 A \right) \EE\left[\| z^{k+1/2} - w^k \|^2\right] + 2\gamma^2 B \EE\left[\sigma^2_k\right] \\
&\hspace{0.4cm} - 2\gamma \EE\left[\langle F(z^{k+1/2}) - F(z^{*}),  z^{k+1/2} - z^*  \rangle \right] + 2\gamma^2 D_1 .
\end{align*}
И по Предположению \ref{as1} в сильно монотонном случае получим
\begin{align*}
(1 &+\gamma  \mu_h ) \EE\left[\| z^{k+1} - z^*\|^2 \right] \leq \tau \EE\left[\| z^k - z^* \|^2\right] + \left(1 -\tau \right) \EE\left[\| w^k - z^* \|^2\right] \\
&\hspace{0.4cm} - \tau \EE\left[\| z^{k+1/2} - z^k\|^2\right]  - \left(\frac{1}{2} + \gamma \mu_h\right) \EE\left[\|z^{k+1} - z^{k+1/2}\|^2\right] \\
&\hspace{0.4cm} - \left((1 -\tau) - 2\gamma^2 A \right) \EE\left[\| z^{k+1/2} - w^k \|^2\right] + 2\gamma^2 B \EE\left[\sigma^2_k\right] \\
&\hspace{0.4cm} - 2\gamma \mu_F \EE\left[\| z^{k+1/2} - z^*  \|^2\right] + 2\gamma^2 D_1.
\end{align*}
С другой стороны
\begin{align*}
\EE\left[\| w^{k+1} - z^*\|^2\right] =  (1 - \tau) \EE\left[\| z^{k+1} - z^*\|^2\right] + \tau \EE\left[\| w^{k} - z^*\|^2\right].
\end{align*}
Суммируя два предыдущих неравенства:
\begin{align*}
\tau &\EE\left[\| z^{k+1} - z^*\|^2 \right] + \EE\left[\| w^{k+1} - z^*\|^2\right] \leq \tau \EE\left[\| z^k - z^* \|^2\right] \\
&\hspace{0.4cm} + \EE\left[\| w^k - z^* \|^2\right] - \tau \EE\left[\| z^{k+1/2} - z^k\|^2\right] - \gamma \mu_h  \EE\left[\| z^{k+1} - z^*\|^2 \right] \\
&\hspace{0.4cm} - \left(\frac{1}{2} + \gamma \mu_h\right)  \EE\left[\|z^{k+1} - z^{k+1/2}\|^2\right] \\
&\hspace{0.4cm} - \left((1 -\tau) - 2\gamma^2 A\right) \EE\left[\| z^{k+1/2} - w^k \|^2\right] \\ 
&\hspace{0.4cm} + 2\gamma^2 B \EE\left[\sigma^2_k\right] + 2\gamma^2 D_1 - 2\gamma \mu_F \EE\left[\| z^{k+1/2} -z^* \|^2\right].
\end{align*}
Добавив $\EE\left[\gamma^2 T  \sigma^2_{k+1}\right]$, мы получим функцию Ляпунова с левой стороны:
\begin{align*}
\EE\left[V_{k+1}\right] &= \tau\EE\left[\| z^{k+1} - z^*\|^2 \right] + \EE\left[\| w^{k+1} - z^*\|^2\right] + \EE\left[\gamma^2 T  \sigma^2_{k+1}\right]\\
&\leq \tau \EE\left[\| z^k - z^* \|^2\right] + \EE\left[\| w^k - z^* \|^2\right] - \tau \EE\left[\| z^{k+1/2} - z^k\|^2\right] \\
&\hspace{0.4cm} - 2\gamma \mu_F \EE\left[\| z^{k+1/2} -z^* \|^2\right] - \gamma \mu_h  \EE\left[\| z^{k+1} - z^*\|^2 \right] \\
&\hspace{0.4cm}  - \left(\frac{1}{2} + \gamma \mu_h\right) \EE\left[\|z^{k+1} - z^{k+1/2}\|^2\right] \\
&\hspace{0.4cm} - \left((1 -\tau) - 2\gamma^2 A \right) \EE\left[\| z^{k+1/2} - w^k \|^2\right] + 2\gamma^2 B \EE\left[\sigma^2_k\right]  \\
&\hspace{0.4cm} + \EE\left[\gamma^2 T  \sigma^2_{k+1}\right] + 2\gamma^2 D_1 .
\end{align*}
С Предположением \ref{as2} для $\sigma_{k+1}$ получаем
\begin{align*}
\EE\left[V_{k+1}\right] &\leq \tau \EE\left[\| z^k - z^* \|^2\right] + \EE\left[\| w^k - z^* \|^2\right] \\ 
&\hspace{0.4cm}+ \left(1 - \rho + \frac{2B}{T}\right) \EE\left[\gamma^2 T  \sigma^2_{k}\right] \\
&\hspace{0.4cm}- \tau \EE\left[\| z^{k+1/2} - z^k\|^2\right] - \gamma \mu_h  \EE\left[\| z^{k+1} - z^*\|^2 \right] \\
&\hspace{0.4cm} - \left(\frac{1}{2} + \gamma \mu_h\right) \EE\left[\|z^{k+1} - z^{k+1/2}\|^2\right]\\
&\hspace{0.4cm}  - \left((1 -\tau) - 2\gamma^2 A - \gamma^2 T C \right) \EE\left[\| z^{k+1/2} - w^k \|^2\right] \\
&\hspace{0.4cm} - 2\gamma\mu_F \EE\left[\| z^{k+1/2} -z^* \|^2\right]  + \gamma^2 (2D_1 + T D_2) .
\end{align*}
Используем $-\| z^{k+1} - z^*\|^2 \leq - \frac{1}{2} \|z^{k+1/2} - z^* \|^2 + \|z^{k+1} - z^{k+1/2} \|^2$, что дает:
\begin{align*}
\EE\left[V_{k+1}\right] 
&\leq \tau \EE\left[\| z^k - z^* \|^2\right] + \EE\left[\| w^k - z^* \|^2\right] \\
&\hspace{0.4cm}  + \left(1 - \rho + \frac{2B}{T}\right) \EE\left[\gamma^2 T  \sigma^2_{k}\right] - \tau \EE\left[\| z^{k+1/2} - z^k\|^2\right] \\ 
&\hspace{0.4cm} - \left((1 -\tau) - 2\gamma^2 A - \gamma^2 T C  \right) \EE\left[\| z^{k+1/2} - w^k \|^2\right] \\ 
&\hspace{0.4cm} - \left(\frac{1}{2} + \gamma \mu_h\right) \EE\left[\|z^{k+1} - z^{k+1/2}\|^2\right] \\ 
&\hspace{0.4cm} - \gamma \left( 2 \mu_F + \frac{\mu_h }{2} \right)  \tau \EE\left[\| z^{k+1/2} - z^*\|^2\right] \\
&\hspace{0.4cm} -\gamma \left( 2 \mu_F + \frac{\mu_h }{2} \right) \cdot(1 - \tau) \EE\left[\| z^{k+1/2} - z^*\|^2\right] \\
&\hspace{0.4cm} + \gamma^2 (2D_1 + T D_2)  .
\end{align*}
Из простых фактов: $\|z^{k+1/2} - z^*\|^2   \geq \frac{1}{2}\|z^k - z^*\|^2 - \|z^{k+1/2} - z^k\|^2$ и $\|z^{k+1/2} - z^*\|^2   \geq \frac{1}{2}\|w^k - z^*\|^2 - \|z^{k+1/2} - w^k\|^2$, следует
\begin{align}
\label{temp1}
&\EE\left[V_{k+1}\right] 
\leq \tau \EE\left[\| z^k - z^* \|^2\right] + \EE\left[\| w^k - z^* \|^2\right] + \left(1 - \rho + \frac{2B}{T}\right) \EE\left[\gamma^2 T  \sigma^2_{k}\right] \notag\\
&\hspace{0.4cm} - \left((1 -\tau) - 2\gamma^2 A - \gamma^2 T C  - \gamma\left(2\mu_F + \frac{\mu_h}{2} \right) \cdot(1 - \tau) \right) \EE\left[\| z^{k+1/2} - w^k \|^2\right]  \notag\\ 
&\hspace{0.4cm}-\left(\frac{1}{2} + \gamma \mu_h\right) \EE\left[\|z^{k+1} - z^{k+1/2}\|^2\right] \notag\\ 
&\hspace{0.4cm} - \left(1 -  \gamma \left( 2 \mu_F + \frac{\mu_h }{2} \right) \right)\tau \EE\left[\| z^{k+1/2} - z^k\|^2\right] \notag\\ 
&\hspace{0.4cm} - \gamma \left( \mu_F + \frac{\mu_h }{4} \right)  \tau \EE\left[\| z^{k} - z^*\|^2\right] \notag\\ 
&\hspace{0.4cm} -\gamma \left( \mu_F + \frac{\mu_h }{4} \right) \cdot(1 - \tau) \EE\left[\| w^{k} - z^*\|^2\right] + \gamma^2 (2D_1 + T D_2) .
\end{align}
Далее работаем с предпоследней строкой \eqref{temp1}:
\begin{align*}
-\gamma&\left(\mu_F + \frac{\mu_h}{4} \right) \tau  \EE\left[\| z^{k} - z^*\|^2\right] -\gamma \left(\mu_F + \frac{\mu_h}{4} \right) \cdot (1 -\tau) \EE\left[\| w^{k} - z^*\|^2\right] \\ &= -\frac{\gamma}{2}\left(\mu_F + \frac{\mu_h}{4} \right) \tau  \EE\left[\| z^{k} - z^*\|^2\right] \\
&\hspace{0.4cm}- \frac{\gamma}{2}\left(\mu_F + \frac{\mu_h}{4} \right) \tau  \EE\left[\| z^{k} - z^*\|^2\right] \\
&\hspace{0.4cm}-\gamma \left(\mu_F + \frac{\mu_h}{4} \right) \cdot (1 -\tau) \EE\left[\| w^{k} - z^*\|^2\right]
\\ &\leq -\frac{\gamma}{2}\left(\mu_F + \frac{\mu_h}{4} \right) \tau  \EE\left[\| z^{k} - z^*\|^2\right] - \frac{\gamma}{4}\left(\mu_F + \frac{\mu_h}{4} \right) \tau  \EE\left[\| w^{k} - z^*\|^2\right] \\
&\hspace{0.4cm} + \frac{\gamma}{2}\left(\mu_F + \frac{\mu_h}{4} \right) \tau  \EE\left[\| z^{k} - w^k\|^2\right] \\
&\hspace{0.4cm} -\gamma \left(\mu_F + \frac{\mu_h}{4} \right) \cdot (1 -\tau) \EE\left[\| w^{k} - z^*\|^2\right]
\\ &\leq -\frac{\gamma}{4}\left(\mu_F + \frac{\mu_h}{4} \right) \tau  \EE\left[\| z^{k} - z^*\|^2\right] - \frac{\gamma}{4}\left(\mu_F + \frac{\mu_h}{4} \right)   \EE\left[\| w^{k} - z^*\|^2\right] \\
&\hspace{0.4cm} + \frac{\gamma}{2}\left(\mu_F + \frac{\mu_h}{4} \right) \tau  \EE\left[\| z^{k} - w^k\|^2\right]
\\ &\leq -\frac{\gamma}{4}\left(\mu_F + \frac{\mu_h}{4} \right) \tau  \EE\left[\| z^{k} - z^*\|^2\right] - \frac{\gamma}{4}\left(\mu_F + \frac{\mu_h}{4} \right)   \EE\left[\| w^{k} - z^*\|^2\right] \\
&\hspace{0.4cm} + \gamma\left(\mu_F + \frac{\mu_h}{4} \right) \tau  \EE\left[\| z^{k+1/2} - z^k\|^2\right] \\
&\hspace{0.4cm}+  \gamma\left(\mu_F + \frac{\mu_h}{4} \right) \tau  \EE\left[\| z^{k+1/2} - w^k\|^2\right].
\end{align*}
Подставив в \eqref{temp1}, получим
\begin{align}
\label{temp302}
\EE\left[V_{k+1}\right] 
&\leq \tau \EE\left[\| z^k - z^* \|^2\right] + \EE\left[\| w^k - z^* \|^2\right]  \notag\\
&\hspace{0.4cm} + \left(1 - \rho + \frac{2B}{T}\right) \EE\left[\gamma^2 T  \sigma^2_{k}\right] \notag\\
&\hspace{0.4cm} - \left((1 -\tau) - 2\gamma^2 A - \gamma^2 T C  - \gamma\left(2\mu_F + \frac{\mu_h}{2} \right) \right) \EE\left[\| z^{k+1/2} - w^k \|^2\right]  \notag\\ 
&\hspace{0.4cm} -\left(\frac{1}{2} + \gamma \mu_h\right) \EE\left[\|z^{k+1} - z^{k+1/2}\|^2\right] \notag\\ 
&\hspace{0.4cm} - \left(1 -  3\gamma \left( \mu_F + \frac{\mu_h }{4} \right) \right)\tau \EE\left[\| z^{k+1/2} - z^k\|^2\right] \notag\\ 
&\hspace{0.4cm} - \frac{\gamma}{4} \left( \mu_F + \frac{\mu_h }{4} \right)  \tau \EE\left[\| z^{k} - z^*\|^2\right] \notag\\ 
&\hspace{0.4cm} -\frac{\gamma}{4} \left( \mu_F + \frac{\mu_h }{4} \right) \EE\left[\| w^{k} - z^*\|^2\right] + \gamma^2 (2D_1 + T D_2).
\end{align}
Остается только выбрать $\gamma \leq \min\left\{ \frac{\sqrt{1 - \tau}}{2\sqrt{2A + TC}}; \frac{1 - \tau}{4\left( \mu_F + \mu_h\right)}\right\}$ и $T \geq \frac{4B}{\rho}$ и получить
\begin{align*}
\EE\left[V_{k+1}\right] 
&\leq \left( 1 - \frac{\gamma}{4} \left( \mu_F + \frac{\mu_h }{4} \right) \right) \left(\tau \EE\left[\| z^k - z^* \|^2\right] + \EE\left[\| w^k - z^* \|^2\right] \right)\\
&\hspace{0.4cm} + \left(1 - \frac{\rho}{2}\right) \EE\left[\gamma^2 T  \sigma^2_{k}\right]  + \gamma^2 (2D_1 + T D_2),
\end{align*}
и в результате: 
\begin{align*}
\EE\left[V_{k+1}\right] 
&\leq \max\left\{\left( 1 -\gamma \cdot \frac{\mu_F + \mu_h}{16} \right) ; \left(1 - \frac{\rho}{2} \right) \right\} \cdot \EE\left[V_{k}\right] \\
&\hspace{0.4cm}+ \gamma^2 (2D_1 + T D_2) .
\end{align*}
Выполнение рекурсивных переходов завершает доказательство.

Далее рассмотрим \textbf{Монотонный/выпуклый случай} ($\mu_h =0$, $\mu_F = 0$). Начнем с \eqref{temp0} с дополнительным обозначением $\text{gap}(z^{k+1/2}, u) =\\
 \langle F(z^{k+1/2}),  z^{k+1/2} - u  \rangle + h(z^{k+1/2}) - h(u)$:
\begin{align*}
2\gamma \cdot \text{gap}(z^{k+1/2}, u) &+ \| z^{k+1} - u\|^2 \leq \tau \| z^k - u \|^2 + (1 -\tau) \| w^k - u \|^2 \notag\\
&- \tau \| z^{k+1/2} - z^k\|^2 - (1 -\tau) \| z^{k+1/2} - w^k \|^2 \notag\\ 
& + 2\gamma^2 \| g^{k+1/2} - g^{k}\|^2 \notag\\ 
& - 2\gamma \langle g^{k+1/2} - F(z^{k+1/2}),  z^{k+1/2} - u  \rangle.
\end{align*}
Добавив к обоим частям $\| w^{k+1} - u\|^2 $ и произведя некоторые перестановки, получим
\begin{align*}
2\gamma \cdot \text{gap}(z^{k+1/2}, u)  &\leq \left[\tau \| z^k - u \|^2 + \| w^k - u \|^2\right]\notag\\
& - \left[\tau \| z^{k+1} - u\|^2 + \| w^{k+1} - u\|^2\right] -\tau \| w^k - u \|^2 - \notag\\
&- (1 - \tau) \| z^{k+1} - u\|^2 + \| w^{k+1} - u\|^2 \notag\\
&- \tau \| z^{k+1/2} - z^k\|^2 - (1 -\tau) \| z^{k+1/2} - w^k \|^2 \notag\\ 
& + 2\gamma^2 \| g^{k+1/2} - g^{k}\|^2\\
& - 2\gamma \langle g^{k+1/2} - F(z^{k+1/2}),  z^{k+1/2} - u  \rangle.
\end{align*}
Просуммируем по $k = 0, \ldots, K - 1$, возьмем максимум от обоих частей по $z \in \mathcal{C}$, далее возьмем математическое ожидание и получим
\begin{align*}
2\gamma \cdot &\EE\left[\max_{u \in \mathcal{C}} \sum\limits_{k=0}^{K-1}\text{gap}(z^{k+1/2}, u)  \right]\leq \max_{u \in \mathcal{C}}\left[\tau \| z^0 - u \|^2 + \| w^0 - u \|^2\right]  \notag\\
&\hspace{0.4cm}+\EE\left[\max_{u \in \mathcal{C}} \sum\limits_{k=0}^{K-1} \left[ -\tau \| w^k - u \|^2 - (1 - \tau) \| z^{k+1} - u\|^2 + \| w^{k+1} - u\|^2\right] \right] \notag\\
&\hspace{0.4cm} - \sum\limits_{k=0}^{K-1} \left[\tau \EE\left[\| z^{k+1/2} - z^k\|^2 \right] + (1 -\tau) \EE\left[\| z^{k+1/2} - w^k \|^2 \right] \right. \\
&\hspace{0.4cm} \left. - 2\gamma^2 \EE\left[\| g^{k+1/2} - g^{k}\|^2 \right] \vphantom{\int_1^2} \right] \notag\\ 
& \hspace{0.4cm}+ 2\gamma \EE\left[\max_{u \in \mathcal{C}} \sum\limits_{k=0}^{K-1} \left[ \langle g^{k+1/2} - F(z^{k+1/2}),  u - z^{k+1/2} \rangle \right] \right].
\end{align*}
Используя Предположение \ref{as2} для $\EE\left[\| g^{k+1/2} - g^{k}\|^2 \right]$, получаем
\begin{align*}
2\gamma \cdot &\EE\left[\max_{u \in \mathcal{C}} \sum\limits_{k=0}^{K-1}\text{gap}(z^{k+1/2}, u)  \right]\leq \max_{u \in \mathcal{C}}\left[\tau \| z^0 - u \|^2 + \| w^0 - u \|^2\right]  \notag\\
&\hspace{0.4cm}+\EE\left[\max_{u \in \mathcal{C}} \sum\limits_{k=0}^{K-1} \left[ -\tau \| w^k - u \|^2 - (1 - \tau) \| z^{k+1} - u\|^2 + \| w^{k+1} - u\|^2\right] \right] \notag\\
&\hspace{0.4cm}- \sum\limits_{k=0}^{K-1} \Big[\tau \EE\left[\| z^{k+1/2} - z^k\|^2 \right] + (1 -\tau) \EE\left[\| z^{k+1/2} - w^k \|^2 \right] \\
&\hspace{0.4cm}- 2\gamma^2 \left(A \EE\left[\| z^{k+1/2} -w^k \|^2\right] + B \EE\left[\sigma^2_k\right] + D_1 \right)\Big]\notag\\ 
& \hspace{0.4cm}+ 2\gamma \EE\left[\max_{u \in \mathcal{C}} \sum\limits_{k=0}^{K-1} \left[ \langle g^{k+1/2} - F(z^{k+1/2}),  u - z^{k+1/2} \rangle \right] \right].
\end{align*}
Добавим и вычтем $\sum\limits_{k=0}^{K-1} \Big[ \gamma^2 T \EE\left[\sigma^2_{k+1}\right] \Big]$  и применим Предположение \ref{as2} для $\sigma_k$, что дает
\begin{align*}
2\gamma \cdot &\EE\left[\max_{u \in \mathcal{C}} \sum\limits_{k=0}^{K-1}\text{gap}(z^{k+1/2}, u)  \right]\leq \max_{u \in \mathcal{C}}\left[\tau \| z^0 - u \|^2 + \| w^0 - u \|^2\right]  \notag\\
&\hspace{0.4cm}+\sum\limits_{k=0}^{K-1} \Big[ \gamma^2 T \EE\left[\sigma^2_{k+1}\right] \Big] - \sum\limits_{k=0}^{K-1} \Big[ \gamma^2 T \EE\left[\sigma^2_{k+1}\right] \Big] \notag\\
&\hspace{0.4cm}+\EE\left[\max_{u \in \mathcal{C}} \sum\limits_{k=0}^{K-1} \left[ -\tau \| w^k - u \|^2 - (1 - \tau) \| z^{k+1} - u\|^2 + \| w^{k+1} - u\|^2\right] \right] \notag\\
&\hspace{0.4cm}- \sum\limits_{k=0}^{K-1} \Big[\tau \EE\left[\| z^{k+1/2} - z^k\|^2 \right] + (1 -\tau -2\gamma^2 A) \EE\left[\| z^{k+1/2} - w^k \|^2 \right] \Big]\notag\\
&\hspace{0.4cm}+ 2\gamma^2 K D_1 + \sum\limits_{k=0}^{K-1} \Big[ 2\gamma^2B \EE\left[\sigma^2_k\right] \Big]\notag\\
&\hspace{0.4cm} + 2\gamma \EE\left[\max_{u \in \mathcal{C}} \sum\limits_{k=0}^{K-1} \left[ \langle g^{k+1/2} - F(z^{k+1/2}),  u - z^{k+1/2}\rangle \right] \right] \\
&\leq \max_{u \in \mathcal{C}}\left[\tau \| z^0 - u \|^2 + \| w^0 - u \|^2\right]  \notag\\
&\hspace{0.4cm}+\sum\limits_{k=0}^{K-1} \Big[ \gamma^2 T \left( (1 - \rho)\EE\left[\sigma^2_k\right]  + C \EE\left[\| z^{k+1/2} -w^k \|^2\right] + D_2 \right) \Big] \notag\\
&\hspace{0.4cm} - \sum\limits_{k=0}^{K-1} \Big[ \gamma^2 T \EE\left[\sigma^2_{k+1}\right] \Big] \notag\\
&\hspace{0.4cm}+\EE\left[\max_{u \in \mathcal{C}} \sum\limits_{k=0}^{K-1} \left[ -\tau \| w^k - u \|^2 - (1 - \tau) \| z^{k+1} - u\|^2 + \| w^{k+1} - u\|^2\right] \right] \notag\\
&\hspace{0.4cm}- \sum\limits_{k=0}^{K-1} \Big[\tau \EE\left[\| z^{k+1/2} - z^k\|^2 \right] + (1 -\tau -2\gamma^2 A) \EE\left[\| z^{k+1/2} - w^k \|^2 \right] \Big]\notag\\
&\hspace{0.4cm}+ 2\gamma^2 K D_1 + \sum\limits_{k=0}^{K-1} \Big[ 2\gamma^2B \EE\left[\sigma^2_k\right] \Big]\notag\\
&\hspace{0.4cm} + 2\gamma \EE\left[\max_{u \in \mathcal{C}} \sum\limits_{k=0}^{K-1} \left[ \langle g^{k+1/2} - F(z^{k+1/2}),  u - z^{k+1/2}  \rangle \right] \right] \notag\\
&= \max_{u \in \mathcal{C}}\left[\tau \| z^0 - u \|^2 + \| w^0 - u \|^2\right]  \notag\\
&\hspace{0.4cm}+\sum\limits_{k=0}^{K-1} \Big[ \gamma^2 T \left( 1 + \frac{2B}{T} -\rho \right)\EE\left[\sigma^2_k\right]  \Big] - \sum\limits_{k=0}^{K-1} \Big[ \gamma^2 T \EE\left[\sigma^2_{k+1}\right] \Big] \notag\\
&\hspace{0.4cm}+\EE\left[\max_{u \in \mathcal{C}} \sum\limits_{k=0}^{K-1} \left[ -\tau \| w^k - u \|^2 - (1 - \tau) \| z^{k+1} - u\|^2 + \| w^{k+1} - u\|^2\right] \right] \notag\\
&\hspace{0.4cm}- \sum\limits_{k=0}^{K-1} \Big[\tau \EE\left[\| z^{k+1/2} - z^k\|^2 \right] \notag\\
&\hspace{0.4cm} + (1 -\tau -\gamma^2 (2A+TC)) \EE\left[\| z^{k+1/2} - w^k \|^2 \right] \Big]\notag\\
&\hspace{0.4cm}+ \gamma^2 K (2D_1 + TD_2) \notag\\
&\hspace{0.4cm} + 2\gamma \EE\left[\max_{u \in \mathcal{C}} \sum\limits_{k=0}^{K-1} \left[ \langle g^{k+1/2} - F(z^{k+1/2}),  u - z^{k+1/2}   \rangle \right] \right].
\end{align*}
С $\gamma \leq \frac{\sqrt{1 -\tau}}{\sqrt{2A + TC}}$ и $T \geq \frac{2B}{\rho}$ получим
\begin{align}
\label{temp4041}
2\gamma \cdot &\EE\left[\max_{u \in \mathcal{C}} \sum\limits_{k=0}^{K-1}\text{gap}(z^{k+1/2}, u)  \right]\leq  \max_{u \in \mathcal{C}}\left[\tau \| z^0 - u \|^2 + \| w^0 - u \|^2\right]  + \gamma^2 T \sigma^2_{0}\notag\\
&\hspace{0.4cm}+\EE\left[\max_{u \in \mathcal{C}} \sum\limits_{k=0}^{K-1} \left[ -\tau \| w^k - u \|^2 - (1 - \tau) \| z^{k+1} - u\|^2 + \| w^{k+1} - u\|^2\right] \right] \notag\\
&\hspace{0.4cm} + 2\gamma \EE\left[\max_{u \in \mathcal{C}} \sum\limits_{k=0}^{K-1} \left[ \langle g^{k+1/2} - F(z^{k+1/2}),  u - z^{k+1/2}  \rangle \right] \right] \notag\\
&\hspace{0.4cm} + \gamma^2 K (2D_1 + TD_2).
\end{align}

Для того, чтобы завершить доказательство, надо оценить члены в последних двух строках. Начнем с $\E\left[ \max\limits_{u \in \mathcal{C}}\sum\limits_{k=0}^{K-1} \langle F(z^{k+1/2}) - g^{k+1/2},  z^{k+1/2} - u \rangle\right]$. Определим последовательность $v$: $v^0 = z^{0}$, $v^{k+1} = \text{prox}_{\gamma h}(v^k-\gamma \delta_k)$ с $\delta^k = F(z^{k+1/2}) - g^{k+1/2}$. Тогда  получим
\begin{align}
    \label{t4441}
    \sum\limits_{k=0}^{K-1} \langle \delta^k, z^{k+1/2} - u \rangle = \sum\limits_{k=0}^{K-1} \langle \delta^k, z^{k+1/2} - v^k \rangle + 
    \sum\limits_{k=0}^{K-1} \langle \delta^k,  v^k - u \rangle . 
\end{align}
По определению $v^{k+1}$ (свойство prox), для всех $z \in \mathcal{Z}$
\begin{align*}
    \langle v^{k+1} - v^{k} +\gamma\delta^k, z - v^{k+1} \rangle \geq 0.
\end{align*}
Переписав этого неравенство,  получим
\begin{align*}
    \langle \gamma\delta^k, v^k  - z \rangle &\leq \langle \gamma\delta^k, v^k - v^{k+1} \rangle  + \langle v^{k+1} - v^k, z - v^{k+1} \rangle \nonumber\\
&\leq \langle \gamma\delta^k, v^k - v^{k+1} \rangle + \frac{1}{2}\|v^k - z\|^2 \\
&\hspace{0.4cm}-  \frac{1}{2}\|v^{k+1} - z\|^2 - \frac{1}{2}\| v^k - v^{k+1}\|^2
\nonumber\\
&\leq \frac{\gamma^2}{2} \|\delta^k\|^2  + \frac{1}{2}\|v^k - v^{k+1}\|^2 + \frac{1}{2}\|v^k - z\|^2 -  \frac{1}{2}\|v^{k+1} - z\|^2
\nonumber\\
&\hspace{0.4cm} - \frac{1}{2}\| v^k - v^{k+1}\|^2 \nonumber\\
&= \frac{\gamma^2}{2} \|\delta^k\|^2  + \frac{1}{2}\|v^k - z\|^2 -  \frac{1}{2}\|v^{k+1} - z\|^2 .
\end{align*}
Вместе с \eqref{t4441} это дает
\begin{align*}
    \sum\limits_{k=0}^{K-1} \langle \delta^k, z^{k+1/2} - u \rangle &\leq \sum\limits_{k=0}^{K-1} \langle \delta^k, z^{k+1/2} - v^k \rangle \\
    & + \frac{1}{\gamma}\sum\limits_{k=0}^{K-1} \left(\frac{\gamma^2}{2} \|\delta^k\|^2  + \frac{1}{2}\|v^k - u\|^2 -  \frac{1}{2}\|v^{k+1} - u\|^2 \right) \nonumber\\
    &\leq \sum\limits_{k=0}^{K-1} \langle \delta^k, z^{k+1/2} - v^k \rangle + 
    \frac{\gamma}{2}\sum\limits_{k=0}^{K-1} \|\delta^k\|^2 + \frac{1}{2\gamma}\|v^0 - u\|^2. 
\end{align*}

Берем максимум по $u$ и получаем
\begin{align*}
     \max_{u \in \mathcal{C}} \sum\limits_{k=0}^{K-1} \langle \delta^k, z^{k+1/2} - u \rangle &\leq \sum\limits_{k=0}^{K-1} \langle \delta^k, z^{k+1/2} - v^k \rangle + \frac{1}{2\gamma} \max_{u \in \mathcal{C}}  \|v^0 - u\|^2 \notag\\
&\hspace{0.4cm} + 
    \frac{\gamma}{2}\sum\limits_{k=0}^{K-1} \|F(z^{k+1/2}) - g^{k+1/2}\|^2 .  
\end{align*}
Берем полное математическое ожидание и получаем
\begin{align}
\label{temp404}
     \EE&\left[ \max_{u \in \mathcal{C}} \sum\limits_{k=0}^{K-1} \langle \delta^k, z^{k+1/2} - u \rangle\right] \leq \EE\left[\sum\limits_{k=0}^{K-1} \langle \delta^k, z^{k+1/2} - v^k \rangle\right] \nonumber\\
    &\hspace{0.4cm}+ 
    \frac{\gamma}{2}\sum\limits_{k=0}^{K-1} \EE\left[\|F(z^{k+1/2}) - g^{k+1/2}\|^2 \right] +\frac{1}{2\gamma}\max_{u \in \mathcal{C}} \|v^0 - u\|^2 \nonumber\\
    &= \EE\left[\sum\limits_{k=0}^{K-1} \langle \EE\left[F(z^{k+1/2}) - g^{k+1/2} ~|~ z^{k+1/2} - v^k \right], z^{k+1/2} - v^k \rangle\right] \nonumber\\
    &\hspace{0.4cm}+ 
    \frac{\gamma}{2}\sum\limits_{k=0}^{K-1} \EE\left[\|F(z^{k+1/2}) - g^{k+1/2}\|^2 \right]+ \frac{1}{2\gamma}\max_{u \in \mathcal{C}} \|v^0 - u\|^2 \nonumber\\
    &= \frac{\gamma}{2}\sum\limits_{k=0}^{K-1} \EE\left[\|F(z^{k+1/2}) - g^{k+1/2}\|^2 \right]+ \frac{1}{2\gamma}\max_{u \in \mathcal{C}} \|v^0 - u\|^2.
\end{align}
Далее оценим $$\EE\left[\max\limits_{u \in \mathcal{C}} \sum\limits_{k=0}^{K-1} \left[ -\tau \| w^k - u \|^2 - (1 - \tau) \| z^{k+1} + u\|^2 + \| w^{k+1} - u\|^2\right] \right],$$ для этого заметим, что
\begin{align*}
    \EE&\left[\max_{u \in \mathcal{C}} \sum\limits_{k=0}^{K-1} \left[ -\tau \| w^k - u \|^2 - (1 - \tau) \| z^{k+1} - u\|^2 + \| w^{k+1} - u\|^2\right] \right] \\
    &= \EE\left[\max_{u \in \mathcal{C}} \sum\limits_{k=0}^{K-1} \left[ -2 \langle (1 - \tau) z^{k+1} + \tau w^k - w^{k+1}, u\rangle \right. \right. \\
    &\hspace{0.4cm} \left. \left.  - (1 -\tau)\|z^{k+1} \|^2 - \tau \| w^k\|^2 + \|w^{k+1} \|^2\right] \vphantom{\int^2} \right] \\
    &= \EE\left[\max_{u \in \mathcal{C}} \sum\limits_{k=0}^{K-1} \left[ -2 \langle (1 - \tau) z^{k+1} + \tau w^k - w^{k+1}, u\rangle \vphantom{\int} \right] \right] \\
    &\hspace{0.4cm} +\EE\left[\sum\limits_{k=0}^{K-1} -(1 -\tau)\|z^{k+1} \|^2 - \tau \| w^k\|^2 + \|w^{k+1} \|^2\right] .
\end{align*}
По определению $w^{k+1}$: $\EE\left[ (1 -\tau)\|z^{k+1} \|^2 + \tau \| w^k\|^2 - \|w^{k+1} \|^2\right] = 0$, тогда 
\begin{align*}
    \EE&\left[\max_{u \in \mathcal{C}} \sum\limits_{k=0}^{K-1} \left[ -\tau \| w^k - u \|^2 - (1 - \tau) \| z^{k+1} - u\|^2 + \| w^{k+1} - u\|^2\right] \right] \\
    &\hspace{4cm}= 2\EE\left[\max_{u \in \mathcal{C}} \sum\limits_{k=0}^{K-1} \langle (1 - \tau) z^{k+1} + \tau w^k - w^{k+1}, -u\rangle\right] 
    \\
    &\hspace{4cm}= 2\EE\left[\max_{u \in \mathcal{C}} \sum\limits_{k=0}^{K-1} \langle (1 - \tau) z^{k+1} + \tau w^k - w^{k+1}, u\rangle\right] .
\end{align*}
Далее можно провести рассуждения аналогично цепочке рассуждений для \eqref{temp404}:
\begin{align}
\label{temp4042}
    \EE&\left[\max_{u \in \mathcal{C}} \sum\limits_{k=0}^{K-1} \left[ \tau \| w^k - u \|^2 + (1 - \tau) \| z^{k+1} - u\|^2 - \| w^{k+1} - u\|^2\right] \right] \notag\\
    &\leq \sum\limits_{k=0}^{K-1} \EE\left[\|(1 - \tau) z^{k+1} + \tau w^k - w^{k+1}\|^2 \right]+ \max_{u \in \mathcal{C}} \|v^0 - u\|^2 \notag\\
    &= \sum\limits_{k=0}^{K-1} \EE\left[\|\EE_{w^{k+1}} [w^{k+1}] - w^{k+1}\|^2 \right]+ \max_{u \in \mathcal{C}} \|v^0 - u\|^2 \notag\\
    &= \sum\limits_{k=0}^{K-1} \EE\left[-\|\EE_{w^{k+1}} [w^{k+1}] \|^2  + \EE_{w^{k+1}} \| w^{k+1}\|^2 \right]+ \max_{u \in \mathcal{C}} \|v^0 - u\|^2 \notag\\
    &= \sum\limits_{k=0}^{K-1} \EE\left[-\|(1 - \tau) z^{k+1} + \tau w^k \|^2  + (1 - \tau)\| z^{k+1}\|^2 + \tau\| w^{k}\|^2 \right] \notag\\
    &+ \max_{u \in \mathcal{C}} \|v^0 - u\|^2 = \sum\limits_{k=0}^{K-1} \tau(1 - \tau)\EE\left[\| z^{k+1} - w^k \|^2 \right]+ \max_{u \in \mathcal{C}} \|v^0 - u\|^2.
\end{align}
Подставив \eqref{temp404} и \eqref{temp4042} в \eqref{temp4041} получим
\begin{align*}
2\gamma \cdot &\EE\left[\max_{u \in \mathcal{C}} \sum\limits_{k=0}^{K-1}\text{gap}(z^{k+1/2}, u)  \right]\leq  \max_{u \in \mathcal{C}}\left[(2 + \tau) \| z^0 - u \|^2 + \| w^0 - u \|^2\right]  \notag\\
&\hspace{0.4cm} + \gamma^2 T \sigma^2_{0} \notag\\
&\hspace{0.4cm} +\sum\limits_{k=0}^{K-1} \left[\tau(1 - \tau)\EE\left[\| z^{k+1} - w^k \|^2 \right] + \gamma^2\EE\left[\|F(z^{k+1/2}) - g^{k+1/2}\|^2 \right] \right]\notag\\
&\hspace{0.4cm}  + \gamma^2 K (2D_1 + TD_2).
\end{align*}
Предположение \ref{as2} для $\EE\left[\|F(z^{k+1/2}) - g^{k+1/2}\|^2 \right]$ дает
\begin{align*}
2\gamma \cdot &\EE\left[\max_{u \in \mathcal{C}} \sum\limits_{k=0}^{K-1}\text{gap}(z^{k+1/2}, u)  \right]\leq  \max_{u \in \mathcal{C}}\left[(2 + \tau) \| z^0 - u \|^2 + \| w^0 - u \|^2\right] \notag\\
&\hspace{0.4cm} +\sum\limits_{k=0}^{K-1} \left[\tau(1 - \tau) \EE\left[\| z^{k+1} - w^k \|^2 \right] + \gamma^2 E \EE\left[\| z^{k+1/2} -w^k \|^2\right] \right] \notag\\
&\hspace{0.4cm} + \gamma^2 T \sigma^2_{0} + \gamma^2 K (2D_1 + TD_2 + D_3).
\end{align*}
С $\gamma \leq \frac{\sqrt{1 - \tau}}{\sqrt{E}}$ приходим к
\begin{align*}
2\gamma \cdot &\EE\left[\max_{u \in \mathcal{C}} \sum\limits_{k=0}^{K-1}\text{gap}(z^{k+1/2}, u)  \right]\leq  \max_{u \in \mathcal{C}}\left[(2 + \tau) \| z^0 - u \|^2 + \| w^0 - u \|^2\right]  \notag\\
&\hspace{0.4cm}  +(1 - \tau) \sum\limits_{k=0}^{K-1} \left[\EE\left[\| z^{k+1} - w^k \|^2 \right] + \EE\left[\| z^{k+1/2} -w^k \|^2\right] \right] \notag\\
&\hspace{0.4cm} + \gamma^2 T \sigma^2_{0} + \gamma^2 K (2D_1 + TD_2 + D_3) \\
&\leq  \max_{u \in \mathcal{C}}\left[(2 + \tau) \| z^0 - u \|^2 + \| w^0 - u \|^2\right]  + \gamma^2 T \sigma^2_{0}\notag\\
&\hspace{0.4cm}+3(1 - \tau) \sum\limits_{k=0}^{K-1} \left[\EE\left[\| z^{k+1} - z^{k+1/2} \|^2 \right] + \EE\left[\| z^{k+1/2} -w^k \|^2\right] \right] \notag\\
&\hspace{0.4cm} + \gamma^2 K (2D_1 + TD_2 + D_3) .
\end{align*}
Вернемся к \eqref{temp302} с $\mu_h = 0$, $\mu_F = 0$, $T \geq \frac{2B}{\rho}$, $\gamma \leq \frac{\sqrt{1-\tau}}{2\sqrt{2A + TC}}$ и получим
\begin{align*}
\EE\left[V_{k+1}\right] 
&\leq \EE\left[V_{k}\right]  - \left((1 -\tau) - 2\gamma^2 A - \gamma^2 T C \right) \EE\left[\| z^{k+1/2} - w^k \|^2\right]  \notag\\ 
&\hspace{0.4cm} -\frac{1}{2} \EE\left[\|z^{k+1} - z^{k+1/2}\|^2\right] + \gamma^2 (2D_1 + T D_2) \\
&\leq \EE\left[V_{k}\right]  - \frac{(1 -\tau)}{2}\EE\left[\| z^{k+1/2} - w^k \|^2\right]  \notag\\ 
&\hspace{0.4cm} -\frac{(1 -\tau)}{2} \EE\left[\|z^{k+1} - z^{k+1/2}\|^2\right] + \gamma^2 (2D_1 + T D_2).
\end{align*}
Следовательно, подставляя это,  подходим к концу доказательства:
\begin{align*}
2\gamma \cdot &\EE\left[\max_{u \in \mathcal{C}} \sum\limits_{k=0}^{K-1}\text{gap}(z^{k+1/2}, u)  \right]
\leq  \max_{u \in \mathcal{C}}\left[(2 + \tau) \| z^0 - u \|^2 + \| w^0 - u \|^2\right] \notag\\
&\hspace{0.4cm} + \gamma^2 T \sigma^2_{0} +6 \sum\limits_{k=0}^{K-1} \left[ \EE\left[V_{k}\right] - \EE\left[V_{k+1} \right] + \gamma^2 (2D_1 + T D_2)\right] \notag\\
&\hspace{0.4cm}  + \gamma^2 K (2D_1 + TD_2 + D_3)\\
&\leq  \max_{u \in \mathcal{C}}\left[(2 + 7\tau) \| z^0 - u \|^2 + 7\| w^0 - u \|^2\right]  + 7\gamma^2 T \sigma^2_{0}  \\
& + \gamma^2 K (14 D_1 + 7TD_2 + D_3)\\
&\leq  \max_{u \in \mathcal{C}}\left[16 \| z^0 - u \|^2\right]  + 7\gamma^2 T \sigma^2_{0}  + \gamma^2 K (14 D_1 + 7TD_2 + D_3).
\end{align*}
 Остается немного подкорректировать критерий сходимости на монотонность $F$ и неравенство Йенсена для выпуклых функций:
\begin{align*}
\EE&\left[\max_{u \in \mathcal{C}} \sum\limits_{k=0}^{K-1}\text{gap}(z^{k+1/2}, u)  \right]\\
& = \EE\left[\max_{u \in \mathcal{C}} \sum\limits_{k=0}^{K-1} \left[ \langle F(z^{k+1/2}),  z^{k+1/2} - u  \rangle + h(z^{k+1/2}) - h(u) \right]  \right] \\
&\geq \EE\left[\max_{u \in \mathcal{C}} \sum\limits_{k=0}^{K-1} \left[ \langle F(u),  z^{k+1/2} - u  \rangle + h(z^{k+1/2}) - h(u) \right]  \right]
\\
&\geq \EE\left[K \cdot \max_{u \in \mathcal{C}} \left[ \langle F(u),  \bar z^{K} - u  \rangle + h(\bar z^{K}) - h(u) \right]  \right]\\
&= K \cdot \EE\left[ \text{Gap} (\bar z^K)\right],
\end{align*}
где мы используем $\bar z^{K} = \frac{1}{K}\sum\limits_{k=0}^{K-1} z^{k+1/2} $. Что приводит к
\begin{align*}
\EE\left[ \text{Gap} (\bar z^K)\right] 
\leq  \frac{8\max_{u \in \mathcal{C}}\left[\| z^0 - u \|^2\right]  + 4\gamma^2 T \sigma^2_{0}}{\gamma K}  + \gamma  (7 D_1 + 3TD_2 + D_3).
\end{align*}

\subsection{Анализ для различных методов}

В этом разделе  устанавливаем связь между  единым анализом и конкретными методами, удовлетворяющими предположению 2. В скобках указаны разделы Приложения, где представлен псевдокод соответствующего метода, а также его анализ при предположении 2, основанный на применении теоремы 1. В большинстве случаев анализируются операторы, удовлетворяющие следующему условию:
\begin{assumption} \label{as3}
$F(z)$~--- \textit{ограниченно-липшицева} с константами $L$ и $D$, то есть для любых $z_1, z_2 \in \mathcal{Z}$ верно
\begin{align}
\label{aslip}
    \|F(z_1) - F(z_2)\|^2 \leq L^2\|z_1-z_2\|^2 + D^2.
\end{align}
\end{assumption}
Заметим, что для $D = 0$ это эквивалентно определнию липшицевости. При $D> 0$, это предположение покрывает случай, когда оператор не является липшицевым, но ограничен.

$\bullet$ \textbf{Существующие методы} (\ref{sec:eg} - \ref{sec:vreg}). Прежде всего хотелось бы упомянуть методы, которые соответствуют нашему параметризованному предположению. Это, конечно, классический экстраградиентный \cite{juditsky2011solving} для задачи \eqref{VI} + \eqref{stoc}. Далее, отметим методы с однократным вызовом оракула \cite{hsieh2019convergence}; отличие этих методов от классического экстраградиентного в том, что на каждой итерации они вычисляют новое значение оператора $F$ лишь один раз. Например, этого можно добиться, используя значение $F$ с предыдущей итерации: $g^k = F(z^{k-1/2}), g^{k+1/2} = F(z^{k+1/2})$ (в экстраградиентном методе имеем: $g^k = F(z^{k}), g^{k+1/2} = F(z^{k+1/2})$). Вариант метода редукции дисперсии \cite{alacaoglu2021stochastic}, специализированный для задач решения ВН \eqref{VI} + \eqref{MK} также удовлетворяет условиям предлагаемого анализа.

$\bullet$ \textbf{Ко-коэрцитивность.} Это предположение аналогично липшицевости оператора:
\begin{align}
    \|F(z_1) - F(z_2)\|^2 \leq l \langle F(z_1) - F(z_2), z_1-z_2 \rangle.
\end{align}
Легко видеть, что $l$-коэрцивитный оператор также является $l$-липшицевым (обратное, вообще говоря, неверно). Более того, если $F$~--- градиент выпуклой функции, то $l$-липшицевость и $l$-коэрцитивность эквивалентны. В литературе имеется анализ некоторых методов (например, метода редукции дисперсии \cite{chavdarova2019reducing}) с этим дополнительным допущением. Довольно легко проанализировать многие методы решения ВН с предположением ко-коэрцитивности. Мы также могли бы построить унифицированную теорию вокруг него и так перенести многие методы минимизации в контекст ВН. Но основная проблема предположения о ко-коэрцитивности состоит в том, что это свойство не выполняется для самой распространенной, билинейной, задачи. Поэтому такой анализ будет справедлив только для минимизации, а это уже проделано в работе \cite{gorbunov2020unified}.

$\bullet$ {\tt Coord-ES} для \eqref{VI} + \eqref{stoc} (\ref{sec:ceg}). Наш первый новый метод позволяет работать не с полным оператором $F$, а выбирать его случайную координату (координаты) и делать шаг только вдоль нее. Методы этого типа называются координатными \cite{wright2015coordinate}. С помощью таких методов можно произвести более тщательный поиск решения~--- выбрать направления, в которых оператор изменяется в большей степени, и проделывать больше шагов в этих направлениях \cite{nesterov2012efficiency}. Также координатный метод очень близок к безградиентным методам \cite{sadiev2020zeroth}, которые актуальны, когда мы работаем с функциями в соответствии с моделью черного ящика, и не можем вычислить оператор $F$/градиент.

$\bullet$ {\tt Quant-ES} для \eqref{VI} + \eqref{stoc} (\ref{sec:qeg}). Суть {\tt Quant-ES} заключается в использовании так называемого оператора квантизации:
\begin{align*}
    \EE{Q(x)} = x,\quad \EE{\| Q(x) \|^2} = \omega \| x\|, \quad \text{для любых} ~~~x. 
\end{align*}
Такие операторы могут быть рандомизированными или детерминированными, с большим или малым параметром $\omega$ \cite{alistarh2017qsgd,beznosikov2020biased}, но все они имеют одну и ту же функцию - сжать вектор $x$. Методы с квантизацией популярны с точки зрения распределенной оптимизации, поскольку основным узким местом там является коммуникация, а сжатие позволяет передавать меньше информации и, следовательно, выигрывать в этом отношении. Мы представляем метод для вариационных неравенств, который может использовать квантизованный оператор.

$\bullet$ {\tt QVR-ES} for \eqref{VI} + \eqref{MK} (\ref{sec:qvreg}). {\tt QVR-ES} сочетает методы редукции дисперсии и квантизации, т.е. сначала мы выбираем случайную функцию с номером $m$ из $M$ вариантом, а затем также квантизуем ее. В простейшем виде это выглядит так: $Q(F_{m}(z))$ -- в нашем методе это делается немного в другом виде, но суть остается той же. Этот метод красочно демонстрирует гибкость нашего подхода и возможность создания различных комбинаций методов с использованием параметризованного предположения 2.

$\bullet$ {\tt IS-ES} \eqref{VI} + \eqref{MK} (\ref{sec:sreg}). В этом случае мы рассматриваем задачу более общую, чем \eqref{VI} + \eqref{MK}. Здесь мы предполагаем, что мы не вызываем функции случайно и равномерно от $1$ до $M$.
Теперь каждый оператор $F_m$ имеет свой вес, в зависимости от которого его можно вызывать чаще или реже.

$\bullet$ {\tt Local-ES} for \eqref{VI} + \eqref{distr} (\ref{sec:loceg}). Этот метод относится к так называемым локальным методам, которые делают ряд локальных обновлений между периодическими коммуникациями.
Наш метод является рандомизированным \cite{hanzely2020federated} и основан на методе из предыдущего абзаца, и также использует технику сэмплирования по важности.

% $\bullet$ {\tt Clip-ES} \eqref{VI} + \eqref{stoc} . Чтобы избежать нежелательных больших отклонений стохастического значения $F$/градиента, принято применять клиппинг \cite{zhang2019adaptive}. Как отмечалось ранее, этот метод не соответствует нашему параметризованному предположению 2, поскольку клиппинг придает оператору некоторое смещение. Этот метод--- просто бонус, потому что он по настоящему эффективен на практике, особенно для GAN.

%===========================================================================================
\section{Заключение}

%===========================================================================================
%\newpage

\printbibliography

% \begin{thebibliography}{99}
% \bibitem{kelly}
% \textit{Kelly F.P., Maulloo A.K., Tan D.K.H. Rate control for communication networks: shadow prices, proportional fairness and stability}~// Journal of the Operational Research Society. 1998. Vol.~49.  \textnumero~3. P.~237--252.
% \end{thebibliography}

%===========================================================================================
%\begin{appendices}
\newpage
\appendix
\renewcommand{\theequation}{П.\arabic{equation}}
\section*{{$\qquad \qquad \qquad \qquad \qquad \qquad \qquad \qquad \qquad$\itПРИЛОЖЕНИЕ} } \label{sec:apCATD}

% \section{Приложение}

\section{Анализ для различных методов}

% Перед тем как приступить к анализу методов, мы вводим следующее

% \begin{assumption} \label{as3}
% $F(z)$ is \textit{ограниченно-липшицева} с константами $L$ и $D$, то есть для всех $z_1, z_2 \in \mathcal{Z}$ верно
% \begin{align}
% \label{aslip}
%     \|F(z_1) - F(z_2)\|^2 \leq L^2\|z_1-z_2\|^2 + D^2.
% \end{align}
% \end{assumption}

% Заметим, что при $D = 0$ это эквивалентно определению липшицевости. При $D> 0$, это предположение покрывает случай, когда оператор ограничен, но не липшицев. 

\subsection{Экстраградиентный метод} \label{sec:eg}

Начнём с простейшего случая в рамках \eqref{VI} + \eqref{stoc} - стохастического с равномерно ограниченным шумом \cite{juditsky2011solving}:

\begin{equation*}
    F(z) = \EE\left[F (z,\xi)\right], \quad \quad \EE \left[ \| F (z,\xi) - F(z)\|^2\right] \leq \sigma^2,
\end{equation*}
где $z$ и $\xi$ независимые. Для него может быть применён

\begin{algorithm} [th]
	\caption{Экстраградиентный метод ({\tt Extra Step})}
	\label{alg_eg}
	\begin{algorithmic}
\STATE
\noindent {\bf Параметры:}  Размер шага $\gamma$, $K$.\\
\noindent {\bf Инициализация:} Выбрать  $z^0\in \mathcal{Z}$.
\FOR {$k=0,1, 2, \ldots, K-1$ }
\STATE Выбрать случайные $\xi^{k}$, $\xi^{k+1/2}$,
\STATE $z^{k+1/2} = \text{prox}_{\gamma h} (z^k - \gamma F(z^{k}, \xi^{k}))$,
\STATE $z^{k+1} = \text{prox}_{\gamma h} (z^k - \gamma F(z^{k+1/2}, \xi^{k+1/2}))$.
\ENDFOR
	\end{algorithmic}
\end{algorithm}

Заметим, что в этом алгоритме $\tau = 0$, и следовательно $w^k = z^k$ для любого $k$. Также мы полагаем $\sigma_k = 0$. Следующая лемма определяет константы и Предположения \ref{as2}:

\begin{lemma}
Предположим, что $F$ ограниченно-липшицева с константами $L$ и $D$ (Предположение \ref{as3}), тогда $g^k$ и $g^{k+1}$ из Алгоритма \ref{alg_eg} удовлетворяют Предположению \ref{as2} с константами $A = 3L^2$, $D_1 = 3D^2 + 6\sigma^2$, $D_3 = \sigma^2$.
\end{lemma}
\textbf{Доказательство:} Легко убедиться, что $g^{k+1/2}$ несмещённо. Далее,
\begin{align*}
     \EE\left[\|g^{k+1/2} - g^k \|^2 \right] &= \EE\left[\|F(z^{k+1/2}, \xi^{k+1/2}) - F(z^{k}, \xi^{k}) \|^2     \right]\\
     &\leq 3\EE\left[\|F(z^{k+1/2}) - F(z^{k}) \|^2 \right]\\
     &\hspace{0.4cm}+ 3\EE\left[\|F(z^{k+1/2}, \xi^{k+1/2}) - F(z^{k+1/2}) \|^2 \right]\\
     &\hspace{0.4cm}+ 3\EE\left[\| F(z^{k}, \xi^{k}) - F(z^{k}) \|^2 \right] \\
     &\leq 3L^2\EE\left[\|z^{k+1/2} - z^{k} \|^2 \right] + 3D^2 + 6\sigma^2,
\end{align*}
и наконец,
\begin{align*}
     \EE\left[\|g^{k+1/2} - F(z^{k+1/2}) \|^2 \right] &= \EE\left[\|F(z^{k+1/2}, \xi^{k+1/2}) - F(z^{k+1/2}) \|^2\right] \leq \sigma^2,
\end{align*}
\begin{corollary}
Предположим, что $F$ ограниченно-липшицева с константами $L$ и $D$. Тогда \textit{Extra Step} 

$\bullet$ в сильно-монотонном случае с $\gamma \leq \min\left\{ \frac{1}{6 L}; \frac{1}{4\left( \mu_F + \mu_h\right)}\right\}$  удовлетворяет
    \begin{align*}
    \EE\left[\|z^{K} - z^*\|^2\right] 
    &\leq \left( 1 -\gamma \cdot \frac{\mu_F + \mu_h}{16}  \right)^{K-1} \cdot \|z^{0} - z^*\|^2 + \frac{96\gamma (D^2 + 2\sigma^2) }{\mu_F + \mu_h},
    \end{align*}
    
$\bullet$ в монотонном случае $\gamma \leq \frac{1}{3L}$ удовлетворяет
    \begin{align*}
    \EE\left[ \text{Gap} (\bar z^K)\right] 
    \leq  \frac{8\max_{u \in \mathcal{C}}\left[\| z^0 - u \|^2\right]}{\gamma K}  + \gamma  (21 D^2 + 43\sigma^2).
    \end{align*}

\textbf{Замечание.} Этот анализ покрывает гладкий случай при $D = 0$. Чтобы получить оценки для негладкого, но ограниченного оператора $F$, достаточно взять $L = 0$ и положить $\frac{1}{L} = +\infty$.

При правильном выборе $\gamma$ (см., например, \cite{stich2019unified}), можно прийти к следующим оценкам скорости сходимости:
$\bullet$ в сильно-монотонном случае
    \begin{align*}
    \EE\left[\|z^{K} - z^*\|^2\right] 
    &= \mathcal{\tilde O} \left( \exp\left( -\frac{(\mu_F + \mu_h)(K-1)}{96 L}  \right)^{} \cdot \|z^{0} - z^*\|^2\right.\\
    &\left.+ \frac{(D^2 + \sigma^2) }{(\mu_F + \mu_h)^2 (K-1)} \right),
    \end{align*}
    
$\bullet$ в монотонном случае
    \begin{align*}
    \EE\left[ \text{Gap} (\bar z^K)\right] 
    = \mathcal{ O} \left(  \frac{ L \max_{u \in \mathcal{C}}\left[\| z^0 - u \|^2\right]}{K}  + \frac{(D + \sigma) \cdot \max_{u \in \mathcal{C}}\left[\| z^0 - u \|\right]}{\sqrt{K}} \right).
    \end{align*}

\end{corollary}

\subsection{Экстраградиентный метод без дополнительного вызова оракула} \label{sec:peg}

Здесь мы также рассматриваем постановку, аналогичную рассматриваемой в предыдущем разделе: \eqref{VI} + \eqref{stoc}. Однако теперь рассматриваем модификацию метода {\tt Extra Step}.

\begin{algorithm} [th]
	\caption{Экстраградиентный метод без дополнительного вызова оракула ({\tt Past-ES})}
	\label{alg_peg}
	\begin{algorithmic}
\STATE
\noindent {\bf Параметры:}  Размер шага $\gamma$, $K$.\\
\noindent {\bf Инициализация:} Выбрать  $z^0\in \mathcal{Z}$.
\FOR {$k=0,1, 2, \ldots, K-1$ }
\STATE Выбрать случайно $\xi_{k+1/2}$,
\STATE $z^{k+1/2} = \text{prox}_{\gamma h} (z^k - \gamma F(z^{k-1/2}, \xi_{k-1/2}))$,
\STATE $z^{k+1} = \text{prox}_{\gamma h} (z^k - \gamma F(z^{k+1/2}, \xi_{k+1/2}))$,
\ENDFOR
	\end{algorithmic}
\end{algorithm}

%%%%%%%%%%%%%%%%%%%%%%%%%%%%%%%%%%%%%%%%%%%%%%%%%
%%%%%%%%%%%%%%%%%%%%%%%%%%%%%%%%%%%%%%%%%%%%%%%%%

\begin{lemma}
Предположим, что $F$ ограниченно-липшицев с константами $L$ и $M$ (Предположение \ref{as3}), тогда $g^k$ и $g^{k+1}$ из Алгоритма \ref{alg_peg} удовлетворяют Предположению \ref{as2} с константами $\rho = \frac{1}{3}, B=3, C=2L^2, D_1 = 6 \sigma^2, D_2 = 4 D^2 + 12 \sigma^2, D_3 = \sigma^2$.
\end{lemma}
\textbf{Доказательство:}
Положим $\displaystyle \sigma_k^2 = \|F(z^{k-1/2}) - F(z^{k+1/2})\|^2$
\begin{align*}
    \EE \left[\sigma_k^2\right] &\leq 2 \E\left[\|F(z^{k}) - F(z^{k+1/2})\|^2\right] + 2 \E\left[\|F(z^{k-1/2}) - F(z^{k})\|^2\right]\\
    &\leq 2 L^2 \E\left[\|z^{k} - z^{k+1/2}\|^2\right] + 2 L^2 \E\left[\|z^{k-1/2} - z^{k}\|^2\right] + 4D^2\\
    &= 2 L^2 \E\left[\|z^{k} - z^{k+1/2}\|^2\right] + 4D^2\\
    &\qquad+ 2 L^2 \E\left[\|z^{k-1} - \gamma F(z^{k-1/2}, \xi_{k}) - z^{k-1} + \gamma F(z^{k-3/2}, \xi_{k-1})\|^2\right]\\
    &= 2 L^2 \E\left[\|z^{k} - z^{k+1/2}\|^2\right]\\
    &\qquad+ 2 L^2 \gamma^2 \E\left[\|F(z^{k-1/2}, \xi_{k}) - F(z^{k-3/2}, \xi_{k-1})\|^2\right] + 4D^2\\
    &\leq 2 L^2 \E\left[\|z^{k} - z^{k+1/2}\|^2\right]\\
    &\qquad+ 6 L^2 \gamma^2 \E\left[\|F(z^{k-1/2}) - F(z^{k-3/2})\|^2\right] + 4D^2 + 12 \sigma^2\\
    &\leq 2 L^2 \E\left[\|z^{k} - z^{k+1/2}\|^2\right] + \frac{2}{3} \E\left[\sigma_{k-1}^2\right] + 4D^2 + 12 \sigma^2,
\end{align*}
если положить $\gamma \leq \frac{1}{3L}$. Следовательно,
\begin{align*}
    \mathbb{E}\left[\|g^{k+1/2} - g^k\|^2\right] &= \mathbb{E}\left[\|F(z^{k+1/2}, \xi_{k+1}) - F(z^{k-1/2}, \xi_{k})\|^2\right] \\
    &\leq 3 \E\left[\sigma_{k}^2\right] + 6 \sigma^2,
\end{align*}
и наконец,
\begin{align*}
     \EE\left[\|g^{k+1/2} - F(z^{k+1/2}) \|^2 \right] &= \EE\left[\|F(z^{k+1/2}, \xi^{k+1/2}) - F(z^{k+1/2}) \|^2\right] \leq \sigma^2.
\end{align*}

\begin{corollary}
Предположим, что $F$ ограниченно-липшицев с константами $L$ и $D$. Тогда \textit{Past-ES} 

$\bullet$ в сильно-монотонном случае с $T = 36$ и $\gamma \leq \min\left\{ \frac{1}{12 L \sqrt{2}}; \frac{1}{4\left( \mu_F + \mu_h\right)}\right\}$  удовлетворяет
    \begin{align*}
    \EE\left[\|z^{K} - z^*\|^2\right] 
    &\leq \left( 1 -\gamma \cdot \frac{\mu_F + \mu_h}{16}  \right)^{K-1} \cdot \|z^{0} - z^*\|^2 \\
    &\hspace{0.4cm}+ \frac{192 \gamma (37 \sigma^2 + 12 D^2) }{\mu_F + \mu_h},
    \end{align*}
    
$\bullet$ в монотонном случае с $T = 18$ и $\gamma \leq \frac{1}{12 L \sqrt{2}}$ удовлетворяет
    \begin{align*}
    \EE\left[ \text{Gap} (\bar z^K)\right] 
    \leq  \frac{8\max_{u \in \mathcal{C}}\left[\| z^0 - u \|^2\right] + 72 \gamma^2 \sigma_0^2}{\gamma K}  + \gamma  (216 D^2 + 691 \sigma^2).
    \end{align*}

\end{corollary}

%%%%%%%%%%%%%%%%%%%%%%%%%%%%%%%%%%%%%%%%%%%%%%%%%
%%%%%%%%%%%%%%%%%%%%%%%%%%%%%%%%%%%%%%%%%%%%%%%%%

Для следующего метода, мы рассматриваем постановку оптимизации оператора вида суммы: \eqref{VI} +\eqref{MK}.

\subsection{Экстраградиентный метод с редукцией дисперсии} \label{sec:vreg}

\begin{algorithm} [th]
	\caption{Экстраградиентный метод с редукцией дисперсии ({\tt VR-ES})}
	\label{alg_vreg}
	\begin{algorithmic}
\STATE
\noindent {\bf Параметры:}  Размер шага $\gamma$, $K$.\\
\noindent {\bf Инициализация:} Выбрать  $z^0 = w^0 \in \mathcal{Z}$.
\FOR {$k=0,1, 2, \ldots, K-1$ }
\STATE $\bar z^k = \tau z^k + (1 - \tau) w^k$
\STATE Выбрать равномерно случайно $m_k \in 1\ldots M$,
\STATE $z^{k+1/2} = \text{prox}_{\gamma h} (\bar z^k - \gamma F(w^k))$,
\STATE $z^{k+1} = \text{prox}_{\gamma h} (\bar z^k - \gamma( F_{m_k}(z^{k+1/2}) - F_{m_k}(w^k) + F(w^k) ))$,
\STATE $w^{k+1} = \begin{cases}
z^{k+1},&  \text{с вероятностью} ~~ 1 - \tau \\
w^k,& \text{с вероятностью} ~~ \tau
\end{cases}$
\ENDFOR
	\end{algorithmic}
\end{algorithm}

Мы полагаем $\sigma_k = 0$. Следующая лемма даёт значения констант для Предположения \ref{as2}:

\begin{lemma}
Предположим, что каждый $F_{m_k}$ и сам $F$ ограниченно-липшицевы с константами $L$ и $D$ (Предположение \ref{as3}), тогда $g^k$ и $g^{k+1}$ из Алгоритма \ref{alg_vreg} удовлетворяют Предположению \ref{as2} с константами $A = L^2$, $D_1 = D^2$, $E = 4L^2$, $D_3 = 4D^2$.
\end{lemma}
\textbf{Доказательство:} Легко убедиться, что $g^{k+1/2}$ несмещённо. Далее,
\begin{align*}
     &\EE\left[\|g^{k+1/2} - g^k \|^2 \right]\\
     &= \EE\left[\| F_{m_k}(z^{k+1/2}) - F_{m_k}(w^k) + F(w^k) - F(w^k) \|^2     \right]\\
    &= \EE\left[\| F_{m_k}(z^{k+1/2}) - F_{m_k}(w^k) \|^2     \right]\\
     &\leq L^2\EE\left[\|z^{k+1/2} - w^{k} \|^2 \right] + D^2,
\end{align*}
и наконец,
\begin{align*}
     &\EE\left[\|g^{k+1/2} - F(z^{k+1/2}) \|^2 \right]\\ 
     &= \EE\left[\|F_{m_k}(z^{k+1/2}) - F_{m_k}(w^k) + F(w^k) - F(z^{k+1/2}) \|^2\right] \\
     & \leq 2 \EE\left[\|F_{m_k}(z^{k+1/2}) - F_{m_k}(w^k) \|^2\right] + 2 \EE\left[\| F(z^{k+1/2}) - F(w^k) \|^2\right] \\
     & \leq 4 L^2 \EE\left[ \| w^k - z^{k+1/2} \|^2 \right]+ 4 D^2,
\end{align*}

\begin{corollary}
Предположим, что каждый $F_{m_k}$ и сам $F$ ограниченно-липшицев с константами $L$ и $D$. Тогда  \textit{VR-ES}

$\bullet$ в сильно монотонном случае с $\gamma \leq \min\left\{ \frac{\sqrt{1 - \tau}}{2 \sqrt{2} L}; \frac{1-\tau}{4\left( \mu_F + \mu_h\right)}\right\}$  удовлетворяет
    \begin{align*}
    &\EE\left[ \tau\| z^{k+1} - z^*\|^2 +  \| w^{k+1} - z^*\|^2 \right]\\
    &\qquad\leq \left( 1 -\gamma \cdot \frac{\mu_F + \mu_h}{16}  \right)^{K-1} \cdot \left(\tau\| z^{0} - z^*\|^2 + \| w^{0} - z^*\|^2 \right) \\
    &\qquad+ \frac{ 32 \gamma D^2 }{\mu_F + \mu_h},
\end{align*}
    
$\bullet$ в монотонном случае с $\gamma \leq \frac{\sqrt{1 - \tau}}{2 \sqrt{6}L}$ удовлетворяет
    \begin{align*}
    \EE\left[ \text{Gap} (\bar z^K)\right] 
    \leq  \frac{8\max_{u \in \mathcal{C}}\left[\| z^0 - u \|^2\right]}{\gamma K}  + 11 \gamma D^2.
    \end{align*}

\end{corollary}

На каждой итераций, мы вычисляем лишь $1$ оператор из $M$. Но в момент обновления $w^k$, необходимо вычислить все $M$ операторов в новой точке $w^k$. На основании этого, мы можем выбрать оптимальное значение для $\tau$ следующим образом:
$$(1-\tau) M \sim \tau \quad \Rightarrow \quad \tau = \frac{M}{M + 1}.$$

\subsection{Покомпонентный экстраградиентный метод} \label{sec:ceg}

Вернёмся назад и снова рассмотрим наиболее общую постановку без конечных сумм: \eqref{VI}.

\begin{algorithm} [th]
	\caption{Покомпонентный экстраградиентный метод ({\tt Coord-ES})}
	\label{alg_ceg}
	\begin{algorithmic}
\STATE
\noindent {\bf Параметры:}  Размер шага $\gamma$, $K$.\\
\noindent {\bf Инициализация:} Выбрать  $z^0 = w^0 \in \mathcal{Z}$.
\FOR {$k=0,1, 2, \ldots, K-1$ }
\STATE $\bar z^k = \tau z^k + (1 - \tau) w^k$
\STATE Выбрать равномерно случайно $i_k \in 1\ldots d$,
\STATE $z^{k+1/2} = \text{prox}_{\gamma h} (\bar z^k - \gamma F(w^k))$,
\STATE $z^{k+1} = \text{prox}_{\gamma h} (\bar z^k - \gamma ( d[F(z^{k+1/2})]_{i_k} e_{i_k} - d[F(w^k)]_{i_k}  e_{i_k} + F(w^k)))$,
\STATE $w^{k+1} = \begin{cases}
z^{k+1},&  \text{с вероятностью} ~~ 1 - \tau \\
w^k,& \text{с вероятностью} ~~ \tau
\end{cases}$
\ENDFOR
	\end{algorithmic}
\end{algorithm}

Положим $\sigma_k = 0$. Следующая лемма даёт значения констант для Предположения \ref{as2}:

\begin{lemma}
Предположим, что $F$ ограниченно-липшицев с константами $L$ and $D$ (Предположение \ref{as3}), тогда $g^k$ и $g^{k+1}$ из Алгоритма \ref{alg_ceg} удовлетворяют Предположению \ref{as2} с константами $A = d L^2$, $D_1 = d D^2$, $E = 2(d+1)L^2$, $D_3 = 2(d+1)D^2$.
\end{lemma}
\textbf{Доказательство:} Легко убедиться, что $g^{k+1/2}$ несмещённо. Далее,
\begin{align*}
     &\EE\left[\|g^{k+1/2} - g^k \|^2 \right]\\
     &= \EE\left[\| d[F(z^{k+1/2})]_{i_k} e_{i_k} - d[F(w^k)]_{i_k}  e_{i_k} + F(w^k) - F(w^k) \|^2     \right]\\
    &= \EE\left[\| d[F(z^{k+1/2}) - F(w^k)]_{i_k}  e_{i_k} \|^2 \right]\\
     &\leq d \EE\left[\|F(z^{k+1/2}) - F(w^k) \|^2 \right] \leq d
     L^2\EE\left[\|z^{k+1/2} - w^{k} \|^2 \right] + d D^2,
\end{align*}
и наконец,
\begin{align*}
     &\EE\left[\|g^{k+1/2} - F(z^{k+1/2}) \|^2 \right]\\ 
     &= \EE\left[\| d[F(z^{k+1/2})]_{i_k} e_{i_k} - d[F(w^k)]_{i_k}  e_{i_k} + F(w^k) - F(z^{k+1/2}) \|^2\right] \\
     & \leq 2 \EE\left[\|d[F(z^{k+1/2}) - F(w^k)]_{i_k}  e_{i_k} \|^2\right] + 2 \EE\left[\| F(z^{k+1/2}) - F(w^k) \|^2\right] \\
     & \leq 2 ( d + 1) L^2 \EE\left[ \| w^k - z^{k+1/2} \|^2 \right]+ 2 ( d + 1) D^2,
\end{align*}

\begin{corollary}
Предположим, что $F$ ограниченно-липшицев с константами $L$ и $D$. Тогда  \textit{Coord-ES}

$\bullet$ в сильно монотонном случае с $\gamma \leq \min\left\{ \frac{\sqrt{1 - \tau}}{2 \sqrt{2 d} L}; \frac{1-\tau}{4\left( \mu_F + \mu_h\right)}\right\}$ удовлетворяет
    \begin{align*}
    &\EE\left[ \tau\| z^{k+1} - z^*\|^2 +  \| w^{k+1} - z^*\|^2 \right]\\ 
    &\qquad\leq \left( 1 -\gamma \cdot \frac{\mu_F + \mu_h}{16}  \right)^{K-1} \cdot \left(\tau\| z^{0} - z^*\|^2 +  \| w^{0} - z^*\|^2 \right) \\
    &\qquad+ \frac{ 32 \gamma  d D^2 }{\mu_F + \mu_h},
    \end{align*}
    
$\bullet$ в монотонном случае с $\gamma \leq \frac{\sqrt{1 - \tau}}{2 L \sqrt{4d + 2}}$ удовлетворяет
    \begin{align*}
    \EE\left[ \text{Gap} (\bar z^K)\right] 
    \leq  \frac{8\max_{u \in \mathcal{C}}\left[\| z^0 - u \|^2\right]}{\gamma K}  + \gamma (9d + 2) D^2.
    \end{align*}
\end{corollary}
Оптимальное значение $\tau = \frac{d}{d+1}$.

\subsection{Квантизованный экстраградиентный метод} \label{sec:qeg}

В этом разделе рассматривается метод, использующий квантизованный оператор.
\begin{algorithm} [th]
	\caption{Квантизованный экстраградиентный метод ({\tt Quant-ES})}
	\label{alg_qeg}
	\begin{algorithmic}
\STATE
\noindent {\bf Параметры:}  Размер шага $\gamma$, $K$.\\
\noindent {\bf Инициализация:} Выбрать  $z^0 = w^0 \in \mathcal{Z}$.
\FOR {$k=0,1, 2, \ldots, K-1$ }
\STATE $\bar z^k = \tau z^k + (1 - \tau) w^k$
\STATE $z^{k+1/2} = \text{prox}_{\gamma h} (\bar z_k - \gamma F(w^k))$,
\STATE $z^{k+1} = \text{prox}_{\gamma h} (\bar z_k - \gamma Q(F(z^{k+1/2}) - F(w^k)) + \gamma F(w^k))$,
\STATE $w^{k+1} = \begin{cases}
z^{k+1},&  \text{с вероятностью} ~~ 1 - \tau \\
w^k,& \text{с вероятностью} ~~ \tau
\end{cases}$
\ENDFOR
	\end{algorithmic}
\end{algorithm}

\begin{definition}(Квантизация). $Q(x)$ называется квантизацией вектора $x \in \R^d$, если 
\begin{align*}
    \EE{Q(x)} = x,\quad \EE{\| Q(x) \|^2} = \omega \| x\|. 
\end{align*}
для некоторого $\omega$.
\end{definition}

\begin{lemma}
Предположим, что $F$ ограниченно-липшицев с константами $L$ и $D$ (Предположение \ref{as3}), тогда $g^k$ и $g^{k+1}$ из Алгоритма \ref{alg_qeg} удовлетворяют Предположению \ref{as2} с константами $A = \omega L^2$, $D_1 = \omega D^2$, $E = 2(\omega + 1) L^2$, $D_3 = 2(\omega + 1)D^2$.
\end{lemma}
\textbf{Доказательство:} Легко убедиться, что $g^{k+1/2}$ несмещённо. Далее,
\begin{align*}
     \EE\left[\|g^{k+1/2} - g^k \|^2 \right] &= \EE\left[\| Q(F(z^{k+1/2}) - F(w^k)) + F(w^k) - F(w^k) \|^2     \right]\\
    &= \EE\left[\| Q(F(z^{k+1/2}) - F(w^k)) \|^2     \right]\\
     &\leq \omega \EE\left[\| F(z^{k+1/2}) - F(w^k) \|^2     \right] \\
     &\leq \omega L^2 \EE\left[ \| z^{k+1/2} - w^k \|^2 \right]+ \omega D^2
\end{align*}
и наконец,
\begin{align*}
     &\EE\left[\|g^{k+1/2} - F(z^{k+1/2}) \|^2 \right]\\
     &= \EE\left[\|Q(F(z^{k+1/2}) - F(w^k)) + F(w^k) - F(z^{k+1/2}) \|^2\right] \\
     & \leq 2\EE\left[\|Q(F(z^{k+1/2}) - F(w^k)) \|^2\right] + 2\EE\left[\|F(w^k) - F(z^{k+1/2}) \|^2\right]\\
     & \leq 2(\omega + 1) L^2 \EE\left[ \| w^k - z^{k+1/2} \|^2 \right]+ 2(\omega + 1)D^2,
\end{align*}

\begin{corollary}
Предположим, что $F$ ограниченно-липшицев с константами $L$ и $D$. Тогда \textit{Quant-ES}

$\bullet$ в сильно монотонном случае с $\gamma \leq \min\left\{ \frac{\sqrt{1 - \tau}}{2 L \sqrt{2\omega}}; \frac{1-\tau}{4\left( \mu_F + \mu_h\right)}\right\}$ удовлетворяет
    \begin{align*}
    &\EE\left[ \tau\| z^{k+1} - z^*\|^2 +  \| w^{k+1} - z^*\|^2 \right]\\
    &\qquad\leq \left( 1 -\gamma \cdot \frac{\mu_F + \mu_h}{16}  \right)^{K-1} \cdot \left(\tau\| z^{0} - z^*\|^2 +  \| w^{0} - z^*\|^2 \right) \\
    &\qquad+ \frac{ 32 \gamma \omega D^2 }{\mu_F + \mu_h},
    \end{align*}
$\bullet$ в монотонном случае с $\gamma \leq \frac{\sqrt{1 - \tau}}{2 L \sqrt{4 \omega + 2}}$ удовлетворяет
    \begin{align*}
    \EE\left[ \text{Gap} (\bar z^K)\right] 
    \leq  \frac{8\max_{u \in \mathcal{C}}\left[\| z^0 - u \|^2\right]}{\gamma K}  + \gamma (9 \omega + 2) D^2 .
    \end{align*}

\end{corollary}

Рассмотрим случай $D = 0$. Квантизация требуется чтобы сжать информацию, при этом $\omega$ выступает здесь как коэффициент сжатия, то есть мы передаём в $\omega$ раз меньше информации чем если бы мы использовали не квантизованный оператор. Однако раз в $1/ (1-\tau)$ итераций (когда мы обновляем $w^k$) необходимо вычислить именно не квантизованный оператор. Основываясь на этом, можно выбрать оптимальное значение для $\tau$ следующим образом:
$$(1-\tau) \sim \tau \cdot \frac{1}{\omega} \quad \Rightarrow \quad \tau = \frac{\omega}{\omega + 1}.$$

\subsection{Квантизованный экстраградиентный метод с редукцией дисперсии} \label{sec:qvreg}

Далее мы совмещаем идеи квантизации и редукции дисперсии для случая задачи \eqref{VI} + \eqref{MK} с оператором вида конечной суммы.

\begin{algorithm} [th]
	\caption{Квантизованный экстраградиентный метод с редукцией дисперсии}
	\label{alg_qvreg}
	\begin{algorithmic}
\STATE
\noindent {\bf Параметры:}  Размер шага $\gamma$, $K$.\\
\noindent {\bf Инициализация:} Выбрать  $z^0 = w^0 \in \mathcal{Z}$.
\FOR {$k=0,1, 2, \ldots, K-1$ }
\STATE Выбрать равномерно случайно $m_k \in 1\ldots M$,
\STATE $\bar z^k = \tau z^k + (1 - \tau) w^k$
\STATE $z^{k+1/2} = \text{prox}_{\gamma h} (\bar z_k - \gamma F(w^k))$,
\STATE $z^{k+1} = \text{prox}_{\gamma h} (\bar z_k - \gamma Q(F_{m_k}(z^{k+1/2}) - F_{m_k}(w^k)) + \gamma F(w^k))$,
\STATE $w^{k+1} = \begin{cases}
z^{k+1},&  \text{с вероятностью} ~~ 1 - \tau \\
w^k,& \text{с вероятностью} ~~ \tau
\end{cases}$
\ENDFOR
	\end{algorithmic}
\end{algorithm}

Положим $\sigma_k = 0$. Следующая лемма даёт значения констант для Предположения \ref{as2}:

\begin{lemma}
Предположим, что каждый $F_{m_k}$ и сам $F$ ограниченно-липшицевы с константами $L$ и $D$ (Предположение \ref{as3}), тогда $g^k$ и $g^{k+1}$ из Алгоритма \ref{alg_qvreg} удовлетворяют Предположению \ref{as2} с $A = \omega L^2$, $D_1 = \omega D^2$, $E = 2(\omega + 1) L^2$, $D_3 = 2(\omega + 1)D^2$.
\end{lemma}
\textbf{Доказательство:} Легко убедиться, что $g^{k+1/2}$ несмещённо. Далее,
\begin{align*}
     \EE\left[\|g^{k+1/2} - g^k \|^2 \right] &= \EE\left[\| Q(F_{m_k}(z^{k+1/2}) - F_{m_k}(w^k)) + F(w^k) - F(w^k) \|^2     \right]\\
    &= \EE\left[\| Q(F_{m_k}(z^{k+1/2}) - F_{m_k}(w^k)) \|^2 \right] \\
    & \leq \omega\EE\left[\|F_{m_k}(z^{k+1/2}) - F_{m_k}(w^k) \|^2\right] \\
    & \leq \omega L^2 \EE\left[ \| w^k - z^{k+1/2} \|^2 \right]+ \omega D^2,
\end{align*}
и наконец,

\begin{align*}
     &\EE\left[\|g^{k+1/2} - F(z^{k+1/2}) \|^2 \right]\\
     &= \EE\left[\|Q(F_{m_k}(z^{k+1/2}) - F_{m_k}(w^k)) + F(w^k) - F(z^{k+1/2}) \|^2\right] \\
     & \leq 2 \EE\left[\|Q(F_{m_k}(z^{k+1/2}) - F_{m_k}(w^k)) \|^2\right] + 2 \EE\left[\| F(w^k) - F(z^{k+1/2}) \|^2\right]\\
     & \leq 2(\omega + 1)L^2 \EE\left[ \| w^k - z^{k+1/2} \|^2 \right]+ 2(\omega + 1 )D^2 ,
\end{align*}

\begin{corollary}
Предположим, что каждый $F_{m_k}$ и сам $F$ ограниченно-липшицевы с константами $L$ и $D$. Тогда \textit{Квантизованный экстраградиентный метод с редукцией дисперсии}

$\bullet$ в сильно монотонном случае с $\gamma \leq \min\left\{ \frac{\sqrt{1 - \tau}}{2 L \sqrt{2\omega}}; \frac{1-\tau}{4\left( \mu_F + \mu_h\right)}\right\}$ удовлетворяет
    \begin{align*}
    &\EE\left[ \tau\| z^{k+1} - z^*\|^2 +  \| w^{k+1} - z^*\|^2 \right]\\
    &\qquad\leq \left( 1 -\gamma \cdot \frac{\mu_F + \mu_h}{16}  \right)^{K-1} \cdot \left(\tau\| z^{0} - z^*\|^2 +  \| w^{0} - z^*\|^2 \right) \\
    &\qquad+ \frac{ 32 \gamma \omega D^2 }{\mu_F + \mu_h},
    \end{align*}
$\bullet$ в монотонном случае с $\gamma \leq \frac{\sqrt{1 - \tau}}{2 L \sqrt{4 \omega + 2}}$ удовлетворяет
    \begin{align*}
    \EE\left[ \text{Gap} (\bar z^K)\right] 
    \leq  \frac{8\max_{u \in \mathcal{C}}\left[\| z^0 - u \|^2\right]}{\gamma K}  + \gamma (9 \omega + 2) D^2 .
    \end{align*}
\end{corollary}
Оптимальное значение $\tau$ здесь то же, что и в предыдущем разделе.

\subsection{Экстраградиентный метод с сэмплированием по важности} \label{sec:sreg}

Здесь мы рассматриваем более общий случай, нежели просто случай конечной суммы. Теперь каждая функция имеет свой вес $p_m$. А именно, рассмотрим дискретную случайную переменную $\eta$ вида
\begin{align*}
    \PP{(\eta = m)} = p_m, \quad \sum\limits_{m=1}^M p_m = 1.
\end{align*}
На каждом шаге, мы обращаемся к $F_{\eta}$. Веса/вероятности $p_m$ могут быть заданы априори, или же мы можем задать их самостоятельно: например, имеет смысл выбрать $p_m = \frac{L_m}{\sum_m L_m}$.

\begin{algorithm} [th]
	\caption{Экстраградиентный метод с сэмплированием по важности}
	\label{alg_sreg}
	\begin{algorithmic}
\STATE
\noindent {\bf Параметры:}  Размер шага $\gamma$, $K$.\\
\noindent {\bf Инициализация:} Выбрать  $z^0 = w^0 \in \mathcal{Z}$.
\FOR {$k=0,1, 2, \ldots, K-1$ }
\STATE Выбрать $\eta_k$,
\STATE $\bar z^k = \tau z^k + (1 - \tau) w^k$
\STATE $z^{k+1/2} = \text{prox}_{\gamma h} (\bar z_k - \gamma F(w^k))$,
\STATE $z^{k+1} = \text{prox}_{\gamma h} (\bar z_k - \gamma \cdot \frac{1}{p_{\eta_k}}\cdot (F_{\eta_k}(z^{k+1/2}) - F_{\eta_k}(w^k)) + \gamma F(w^k))$,
\STATE $w^{k+1} = \begin{cases}
z^{k+1},&  \text{с вероятностью} ~~ 1 - \tau \\
w^k,& \text{с вероятностью} ~~ \tau
\end{cases}$
\ENDFOR
	\end{algorithmic}
\end{algorithm}

\begin{lemma}
Предположим, что каждый $F_{m}$ ограниченно-липшицев с константами $L_m$ и $D_m$ (Предположение \ref{as3}), как и $F$, соответственно с константами $L$ и $D$, тогда $g^k$ и $g^{k+1}$ из Алгоритма \ref{alg_sreg} удовлетворяют Предположению \ref{as2} с константами $A = \sum\limits_{m=1}^M \frac{L_m^2}{p_m}$, $D_1 = \sum\limits_{m=1}^M \frac{D_m^2}{p_m}$, $E = 2\left( \sum\limits_{m=1}^M \frac{L_m^2}{p_m} + L^2\right)$, $D_3 = 2\left(\sum\limits_{m=1}^M \frac{D_m^2}{p_m} + D^2\right)$.
\end{lemma}
\textbf{Доказательство:} Легко убедиться, что $g^{k+1/2}$ несмещённо. Далее,
\begin{align*}
    &\EE\left[\|g^{k+1/2} - g^k \|^2 \right]\\
    &\qquad= \EE\left[\left\| \frac{1}{p_{\eta_k}} (F_{\eta_k}(z^{k+1/2}) - F_{\eta_k}(w^k)) + F(w^k) - F(w^k) \right\|^2     \right]\\
    &\qquad= \EE\left[\left\| \frac{1}{p_{\eta_k}} (F_{\eta_k}(z^{k+1/2}) - F_{\eta_k}(w^k))\right\|^2     \right]\\
    &\qquad= \EE \sum\limits_{m=1}^M \frac{1}{p_m}\left[\|F_{m}(z^{k+1/2}) - F_{m}(w^k) \|^2\right] \\
    &\qquad \leq \sum\limits_{m=1}^M \frac{L_m^2}{p_m} \EE \left[\|z^{k+1/2} - w^k \|^2\right] + \sum\limits_{m=1}^M \frac{D_m^2}{p_m},
\end{align*}
и наконец,

\begin{align*}
     &\EE\left[\|g^{k+1/2} - F(z^{k+1/2}) \|^2 \right]\\
     &= \EE\left[\left\|\frac{1}{p_{\eta_k}} (F_{m_k}(z^{k+1/2}) - F_{m_k}(w^k)) + F(w^k) - F(z^{k+1/2}) \right\|^2\right] \\
     & \leq 2 \EE\left[\left\|\frac{1}{p_{\eta_k}}(F_{m_k}(z^{k+1/2}) - F_{m_k}(w^k)) \right\|^2\right] + 2 \EE\left[\| F(w^k) - F(z^{k+1/2}) \|^2\right]\\
     & \leq 2\left( \sum\limits_{m=1}^M \frac{L_m^2}{p_m} + L^2\right) \EE\left[ \| w^k - z^{k+1/2} \|^2 \right]+ 2\left(\sum\limits_{m=1}^M \frac{D_m^2}{p_m} + D^2\right) ,
\end{align*}

\begin{corollary}
Предположим, что каждый $F_{m}$ ограниченно-липшицев с константами $L_m$ и $D_m$ (Предположение \ref{as3}) и $F$ с константами $L$ и $D$. Тогда \textit{Экстраградиентный метод с сэмплированием по важности}

$\bullet$ в сильно монотонном случае с $\gamma \leq \min\left\{ \frac{\sqrt{1 - \tau}}{2 \sqrt{2} \sum\limits_{m=1}^M \frac{L_m^2}{p_m} }; \frac{1-\tau}{4\left( \mu_F + \mu_h\right)}\right\}$ удовлетворяет
    \begin{align*}
    &\EE\left[ \tau\| z^{k+1} - z^*\|^2 +  \| w^{k+1} - z^*\|^2 \right]\\
    &\qquad\leq \left( 1 -\gamma \cdot \frac{\mu_F + \mu_h}{16}  \right)^{K-1} \cdot \left(\tau\| z^{0} - z^*\|^2 +  \| w^{0} - z^*\|^2 \right)\\
    &\qquad+ \frac{ 32 \gamma}{\mu_F + \mu_h} \cdot \sum\limits_{m=1}^M \frac{D_m^2}{p_m},
\end{align*}
$\bullet$ в монотонном случае с $\gamma \leq \frac{\sqrt{1 - \tau}}{2 \sqrt{2 L^2 + 4 \sum\limits_{m=1}^M \frac{L_m^2}{p_m} }}$ удовлетворяет
    \begin{align*}
    \EE\left[ \text{Gap} (\bar z^K)\right] 
    \leq  \frac{8\max_{u \in \mathcal{C}}\left[\| z^0 - u \|^2\right]}{\gamma K}  + \gamma \left(9 \sum\limits_{m=1}^M \frac{D_m^2}{p_m} + 2 D^2\right).
    \end{align*}

\end{corollary}

\subsection{Локальный экстраградиентный метод} \label{sec:loceg}

Этот метод предназначен для распределённой задачи \eqref{VI}+\eqref{distr}. Суть метода заключается в переключении между локальными итерациями и усреднением со значением на сервере. С вероятностью $\tau$ производится локальный шаг, с вероятностью $1- \tau$ -- шаг коммуникации. 

Здесь $Z_{\text{avg}} = [\bar z^T, \ldots, \bar z^T] \in \R^{Md}$ с $ \bar z = \frac{1}{M} \sum_{m=1}^M  z_m$ (и то же для $W_{\text{avg}}$).

\begin{algorithm} [H]
\label{alg_loceg}
	\caption{Рандомизированный локальный экстраградиентный метод}
	\begin{algorithmic}
	\STATE
\noindent {\bf Параметры:}  Размер шага $\gamma$, $K$,  вероятность $p$.\\
\noindent {\bf Инициализация:} Выбрать  $z^0 = w^0 \in \mathcal{Z}$, $z^0_m = z^0$ для всех $m$.
\FOR {$k=0,1, 2, \ldots$ }
\STATE $\bar Z^k = \tau Z^k + (1 - \tau) W^k$, 
\STATE $Z^{k+1/2} = \text{prox}_{\gamma h} \left(\bar Z^k - \gamma \cdot (\Phi(W^k) + \lambda  \cdot (W^k - W^k_{\text{avg}})\right)$, 
\STATE $G (Z) = \begin{cases}
\frac{1}{\tau} \cdot \Phi(Z),&  \text{с вероятностью} ~~ \tau \\
\frac{1}{1 - \tau } \cdot \lambda  \cdot (Z - Z_{\text{avg}}) ,& \text{с вероятностью} ~~ 1 - \tau
\end{cases},$
\STATE $Z^{k+1} = \text{prox}_{\gamma h} \left(\bar Z^k - \gamma \cdot (G (Z^{k+1/2}) - G (W^k) + \Phi(W^k) + \lambda  \cdot (W^k - W^k_{\text{avg}}) )\right)$, 
\STATE $W^{k+1} = \begin{cases}
Z^{k+1},&  \text{с вероятностью} ~~ 1 - \tau \\
W^k,& \text{с вероятностью} ~~ \tau
\end{cases}$
\ENDFOR
\end{algorithmic}
\end{algorithm}

На самом деле, анализ для этого метода уже был проделан в предыдущем разделе. Действительно, здесь мы имеем два оператора: $\Phi(Z)$ и $\lambda  \cdot (Z - Z_{\text{avg}})$, являющихся $L$ и $\lambda$ гладкими, соответственно.

Если выбрать $\tau = \frac{L}{L + \lambda}$, то верно

\begin{corollary}
\textit{Рандомизированный локальный экстраградиентный метод}

$\bullet$ в сильно монотонном случае с $\gamma \leq \min\left\{ \frac{\sqrt{\lambda}}{2 \sqrt{2} (L + \lambda)^{3/2}}; \frac{\sqrt{\lambda}}{4\left( \mu_F + \mu_h\right)\sqrt{L + \lambda}}\right\}$  удовлетворяет
    \begin{align*}
    &\EE\left[ \tau\| z^{k+1} - z^*\|^2 + \| w^{k+1} - z^*\|^2 \right]\\
    &\qquad\leq \left( 1 -\gamma \cdot \frac{\mu_F + \mu_h}{16}  \right)^{K-1} \cdot \left(\tau\| z^{0} - z^*\|^2 +  \| w^{0} - z^*\|^2 \right)\\
    &\qquad+ \frac{ 32 \gamma}{\mu_F + \mu_h} \cdot \sum\limits_{m=1}^M \frac{D_m^2}{p_m},
    \end{align*}
$\bullet$ в монотонном случае с $\gamma \leq \frac{\sqrt{\lambda}}{2 \sqrt{6} (L + \lambda)^{3/2}}$ удовлетворяет
    \begin{align*}
    \EE\left[ \text{Gap} (\bar z^K)\right] 
    \leq  \frac{8\max_{u \in \mathcal{C}}\left[\| z^0 - u \|^2\right]}{\gamma K}  + \gamma \left(9 \sum\limits_{m=1}^M \frac{D_m^2}{p_m} + 2 D^2\right).
    \end{align*}
\end{corollary}
Отсюда, можно получить оценку для числа локальных шагов и числа коммуникаций:

$\bullet$ в сильно монотонном случае:
\begin{align*}
\mathcal{ O} \left(\frac{\sqrt{\lambda (\lambda + L)}}{\mu}\log\frac{1}{\varepsilon}\right) \text{комм.~~~и~~~ } \mathcal{ O} \left(\frac{\sqrt{L(\lambda + L)}}{\mu} \log\frac{1}{\varepsilon}\right) \text{локальных шагов. }
 \end{align*}
$\bullet$ в монотонном случае
\begin{align*}
\mathcal{ O} \left(\frac{\sqrt{\lambda (\lambda + L)} \cdot \Omega^2}{\varepsilon}\right) \text{комм.~~~и~~~ } \mathcal{ O} \left(\frac{\sqrt{L (\lambda + L)}\cdot \Omega^2}{\varepsilon} \right) \text{локальных шагов. }
  \end{align*}
Отметим, что эти оценки являются достаточно хорошими при малом значении $\lambda$.

\section{Эксперименты}

\subsection{Генеративно-состязательные сети}

% The purpose of our experiments is to confirm the workability of SVRG, as well as quantization234methods for training GANs.
В качестве примера оптимизации мини-максного целевого функционала при различных подходах оптимизации была предложена задача оптимизации генеративно-состязательных сетей(GAN).  GAN — это структура для оценки генеративных моделей с помощью состязательного процесса, в котором мы одновременно обучаем две модели: генеративную модель $G$, которая фиксирует распределение данных, и дискриминативную модель $D$, которая оценивает вероятность того, что пример пришел из обучающей выборки, а не из $G$. 
${D(G(z))}$ это вероятность (скалярная), что выход генератора $G$ это реальное изображение. Согласно статье Гудфеллоу \cite{goodfellow2014generative}, в этой задаче ${D}$ пытается максимизировать вероятность того, что он правильно классифицирует реальные и фейковые изображения ${(\log(D(x))}$, а ${G}$ пытается минимизировать вероятность того, что ${D} $ будет предсказывать его выходные данные как фейковые $ {(\log(1-D(G(z)))}$. Из статьи функция потерь GAN

\[\min_{\theta_G}\max_{\theta_D} \mathbb{E}_{x\backsim \mathbb{P}_{data}(x)} \big [ \log(D(x)) \big ] + \mathbb{E}_{z\backsim \mathbb{P}_{z}(z)} \big [ \log(1-D(z)) \big ]           \]

Целью наших экспериментов не было получение наилучших результатов для генеративных сетей или подходом к формулировке задачи для GAN. В наших экспериментах мы подтверждаем работоспособность SVRG, а также методов квантизации для обучения GAN.

\textbf{Данные, модель и и оптимизаторы.} 
% We consider the CIFAR-10 dataset. 
Для проведения экспериментов было предложено использовать датасет CIFAR. Он содержит $50000$ и $10000$ изображений в обучающей и валидационной выборках соответственно, равномерно распределенных по $10$ классам.
Для каждого оптимизационного подхода в качестве оптимизатора использовался Adam.
В ходе проведения эксперимента, был проведен анализ гиперпараметроов каждого из оптимизаторов(генератора и дискриминатора):
 \begin{table}[H]
     \centering
     \begin{tabular}{|c|c|c|c|}
     \hline
     $\beta_1$&.5 & .9&.99\\\hline
     $\beta_2$&.99&.999&\\\hline
     \text{скорость обучения}&\ $.9\cdot 10^{-4}$&$1\cdot 10^{-4}$&$2\cdot 10^{-4}$\\\hline
     \end{tabular}
     \caption{Набор анализируемых гиперпараметров}
     \label{gridsearch}
 \end{table}
 По результатам проведенного анализа сходимости (проверялись результаты работы после 100 эпох обучения), были выбраны следующие параметры: $\beta_1=.9\, , \beta_2=.999 \text{ и коэффициент скорости обучения}=2\cdot 10^{-4}$ для обоих методов оптимизации. Размер набора элементов, на котором производилось одна итерация вычисления градиента, составлял 64.
 
 Для проверки качества выдаваемых изображений, была использована метрика Inception score \cite{barratt2018note}.
%  \begin{center}
%  \centering
  \[\textbf{IS}(G)=exp(\mathbb{E}_{x\sim p_g} D_{KL}(p(y|\textbf{x})||p(y)))\]
%  \end{center}
 где $x\sim p_g$ означает что изображение \textbf{x} сгенерировано из распределения $p_g$, $D_{KL}(p||q)$ -- расстояние Кульбака-Лейбнера между двумя распределениями $p,q$.
 
 Данная метрика позволяет автоматически оценивать качество картинок, генерируемых моделью. в ходе исследования было показано, что данная метрика хорошо кореллирует с результатами оценки человека на сгенерированных изображениях CIFAR10. данная метрика использует нейронную сеть Inception v3, предобученную на датасете ImageNet, и собирает статистику выходных данных работы сети, примененной к сгенериорванным изображениям.
% Мы используем Адама в качестве оптимизатора со скоростью обучения $10^{-5}$ для генератора и $10^{-4}$ для дискриминатора.
Для экспериментов было предложено использовать архитектуру DCGan \cite{radford2016unsupervised} с подходом Conditional \cite{mirza2014conditional}. Особенность данной архитектуры заключается в том что мы можем обучить модель сэмплировать изображения из определенного распределения, сэмплы из которого будут схожи с элементами из желаемого распределения обучающей выборки.

Основной целью этого эксперимента была возможность не переобучать определенные архитектуры под разные подходы оптимизации. В процессе обучения генератор и дискриминатор имели одинаковое количество шагов оптимизации.

\begin{figure}
    \centering
    \includegraphics[width =.9\textwidth]{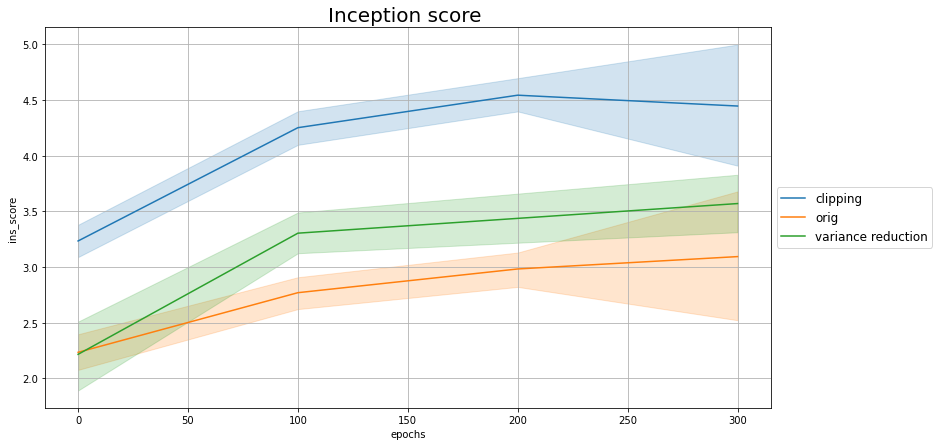}
    \caption{Результаты работы модели, полученных при ее потимизации с помощью метода Adam разными способами: подход с использованием квантизации/клиппинга 70\%  всех градиентов модели, подход с использованием редукции дисперсии, оригинальный подход \cite{goodfellow2014generative}}
    \label{fig:gans1}
\end{figure}

\begin{figure}[h!]
\begin{minipage}{0.5\textwidth}
  \centering
\includegraphics[width =  0.9\textwidth ]{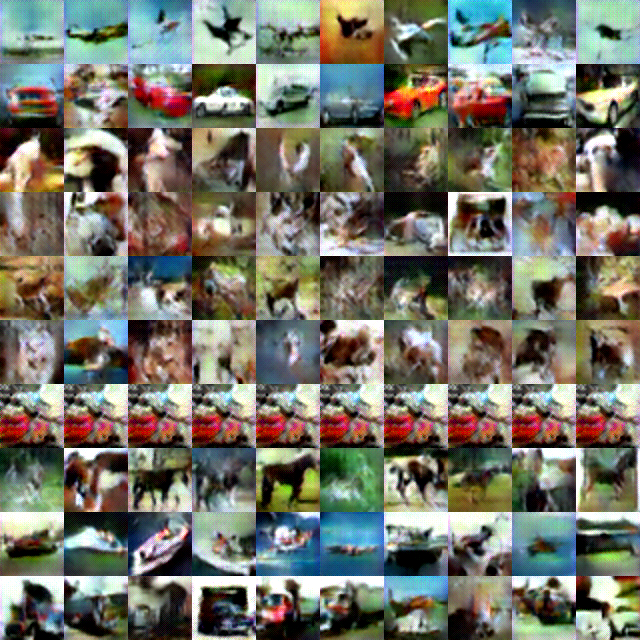}
\end{minipage}%
\begin{minipage}{0.5\textwidth}
  \centering
\includegraphics[width =  0.9\textwidth ]{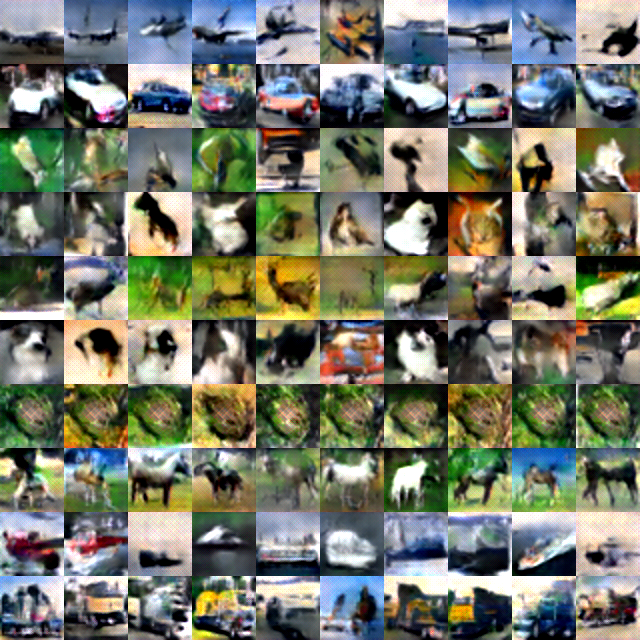}
\end{minipage}%
\\
\begin{minipage}{0.5\textwidth}
  \centering
  (a) 
\end{minipage}%
\begin{minipage}{0.5\textwidth}
  \centering
  (b)
\end{minipage}%
\\
\centering
\begin{minipage}{0.5\textwidth}
  \centering
    \includegraphics[width =  0.9\textwidth]{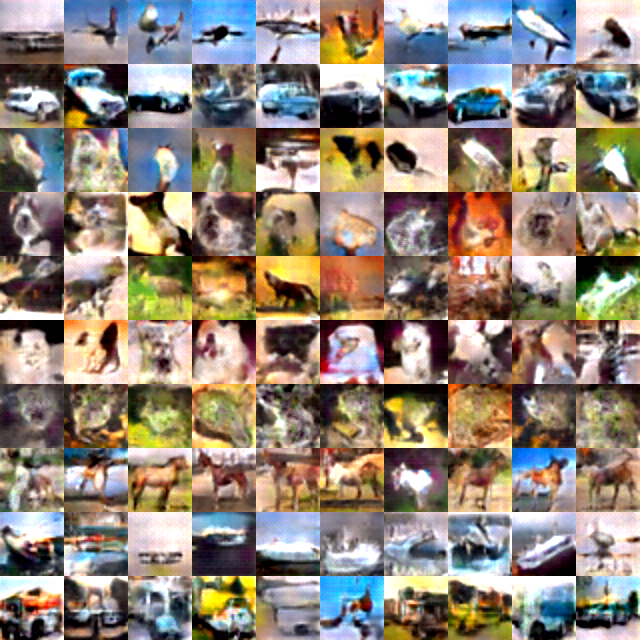}
\end{minipage}%
% % \begin{minipage}{0.5\textwidth}
% %   \centering
% % \includegraphics[width =  0.73\textwidth ]{include/plots/SVRG.png}
% % \end{minipage}%
\\
\centering
\begin{minipage}{0.5\textwidth}
  \centering
  (c)
  
\end{minipage}%

% \begin{minipage}{0.5\textwidth}
%   \centering
%   (d) 
% \end{minipage}%
\caption{Эти изображения были сгенерированы архитектурой conditional DCGan, которая обучалась с использованием разных подходов:(a)-исходная оптимизация с использованием Adam  1, (b) оптимизация с занулением  $70 \%$ всех значений градиентов,
(c) оптимизация со стохастической редукцией дисперсии  градиентного спуска.}
\label{fig:gans2}
\end{figure}

\textbf{Результаты.} Посмотрите на рисунки представленные на фигурах \ref{fig:gans1} и \ref{fig:gans2}.
Согласно результатам, полученным в рамках экспериментов, можно заметить, что подходы, предложенные в данной статье, позволяют достигнуть лучших результатов, по сравнению с оригинальными.

$\bullet$ Метод редукции дисперсии позволял оптимизировать функционал точнее, чем оригинальный способ, на протяжении всей эпохи оптимизации,однако, приближаясь к конечным итерациям каждой эпохи обучения, полный градиент достаточно сильно отличался от направления градиента, получаемого на последних итерациях, что вносило неточность при оптимизации, что привело к переобучению генератора.  

    % overfitting.
$\bullet$  Метод с использованием квантизации/клиппинга (случайным образом обнуляется 70\% от всего градиента модели как в генераторе так и в дискриминаторе, что соответствует оператору сжатия Rand70\%),позволил предотвратить переобучение генератора и дискриминатора на ранних стадиях, что в свою очередь позволило достичь лучших результатов.

% $\bullet$  квантизация всех градиентов позволила предотвратить переобучение дискриминатора и генератора.

% \begin{figure}[h!]
% \begin{minipage}{1.1\textwidth}
%   \centering
% \includegraphics[width =  0.9\textwidth ]{include/plots/gan_experiments.pdf}
% \end{minipage}%
% % \begin{minipage}{0.5\textwidth}
% %   \centering
% % \includegraphics[width =  0.9\textwidth ]{NIPS/local_FID_it.pdf}
% % \end{minipage}%
% % \begin{minipage}{0.5\textwidth}
% %   \centering
% % \end{minipage}%
% % \begin{minipage}{0.5\textwidth}
% %   \centering
% % \end{minipage}%
% \caption{ Сравнение скорости сходимости для разных подходов. Мы сравнили исходную оптимизацию архитектуры DCGAN с использованием Adam с другими подходами к оптимизации, такими как метод редукции дисперсии и зануление (квантизация) большинства значений градиента.}
% \label{fig:gans1}

% \end{figure}

% \begin{figure}[h!]
% \begin{minipage}{0.5\textwidth}
%   \centering
% \includegraphics[width =  0.9\textwidth ]{include/plots/orig_gan.pdf}
% \end{minipage}%
% \begin{minipage}{0.5\textwidth}
%   \centering
% \includegraphics[width =  0.73\textwidth ]{include/plots/cut_descriminator.png}
% \end{minipage}%
% \\
% \begin{minipage}{0.5\textwidth}
%   \centering
%   (a) 
% \end{minipage}%
% \begin{minipage}{0.5\textwidth}
%   \centering
%   (b)
% \end{minipage}%
% \\
% \begin{minipage}{0.5\textwidth}
%   \centering
%     \includegraphics[width =  0.9\textwidth]{include/plots/all_cut_30.pdf}
% \end{minipage}%
% \begin{minipage}{0.5\textwidth}
%   \centering
% \includegraphics[width =  0.73\textwidth ]{include/plots/SVRG.png}
% \end{minipage}%
% \\
% \begin{minipage}{0.5\textwidth}
%   \centering
%   (c)
% \end{minipage}%
% \begin{minipage}{0.5\textwidth}
%   \centering
%   (d) 
% \end{minipage}%
% \caption{Эти изображения были сгенерированы DCGAN, которые обучались с использованием разных подходов:(a)-исходная оптимизация с использованием Adam  1, (b) оптимизация с занулением $70 \%$ значений градиента дискриминатора, (c) оптимизация с занулением  $70 \%$ всех значений градиентов,
% (d) оптимизация со стохастической редукцией дисперсии  градиентного спуска.}
% \label{fig:gans2}
% \end{figure}

\textbf{Подводя итоги}, все вышеупомянутые подходы позволили получить результаты лучше  полученных при оригинальном подходе обучения генеративных моделей и могут быть применены для оптимизации GAN, но квантизация позволяет быстро получить наилучшие результаты при сохранении вычислительных ресурсов.

\subsection{Полицейский vs. Грабитель}

Рассмотрим задачу Полицейский vs. Грабитель, 
% \cite{nemirovski2013mini} 
чтобы сравнить производительность некоторых методов, представленных в нашей статье.
В этой задаче мы рассматриваем квадратный город размером 200 на 200 клеток. В каждом городе есть дом и будка милиции. Грабитель выбирает один дом для ограбления, полицейский выбирает будку, в которой будет дежурить. Задача заключается в нахождении оптимальной смешанной стратегии для противоборствующих игроков: грабителя и полицейского в игре, где цель грабителя -- атаковать выбранный дом $ i $ с максимальным благосостоянием $ w_i $, а цель полицейского -- выбрать оптимальный пост и поймать грабителя, чтобы предотвратить нанесенный им максимальный ожидаемый убыток. Мы предполагаем, что вероятность поймать вора для определенного дома $i$ и поста $j$ равна $\exp\left(-\theta \cdot d(i, j)\right)$ для функции расстояния $d$, которая вводится выше. Мы можем сформулировать описанную постановку как задачу о поиски билинейной седловой точки:
\begin{equation}\label{problem}
    \min_{x \in \Delta(n^2)} \max_{y \in \Delta(n^2)} f(x, y) := \frac{1}{n}\sum_{k=1}^n y^\top A^{(k)} x,
\end{equation}
для $x$ и $y$, которые являются векторами вероятности выбора какого-либо дома и поста соответственно, и для матриц
\[
    A^{(k)}_{i j} = w^{(k)}_i \cdot \left(1 - \exp\left(-\theta \cdot d(i, j)\right)\right),
\]
где благосостояние $w^{(k)}$ и функция расстояния $d$ определяются следующим образом (эти выражения легко понять, если представить $ i $ как сплющенную координату на игровом поле размера $n\times n$, $i(x, y) = x\cdot n + y$, график $w$ представляет собой пирамиду с центром в центре этого поля и $d$ евклидово расстояние на ней):
\begin{align*}
    w_i = 1 - \frac{2}{n} \cdot \min \left\{\left|\left\lfloor i / n\right\rfloor - n/2\right|, \left|i\text{ mod }n - n/2\right|\right\}, \\
    w^{(k)}_i = w_i \cdot (1 + \xi^{(k)})\text{ for }\xi^{(k)} \sim \mathcal{U}(0, \sigma), \\
    d(i, j) = \sqrt{(\left\lfloor i / n\right\rfloor - \left\lfloor j / n\right\rfloor)^2 + (i\text{ mod }n - j\text{ mod }n)^2}.
\end{align*}

Для экспериментов были выбраны параметры $\theta = 0.6, n = 25, \sigma = 3$. {\tt Coord-ES}, {\tt Quant-ES} и {\tt Past-ES} методы использует такой же значение $\gamma$, как {\tt Extra Step} метод, которое является оптимальным для последнего. Мы сравниваем {\tt Past} метод по количеству обращений к оракулу $F$, метод с квантизацией по количеству используемых бит, координатный метод по количеству используемых координат.
% \cite{held1974valid}. 

\begin{figure}[H]
     \centering
     \begin{subfigure}[H]{0.32\textwidth}
         \centering
         \includegraphics[width=\textwidth]{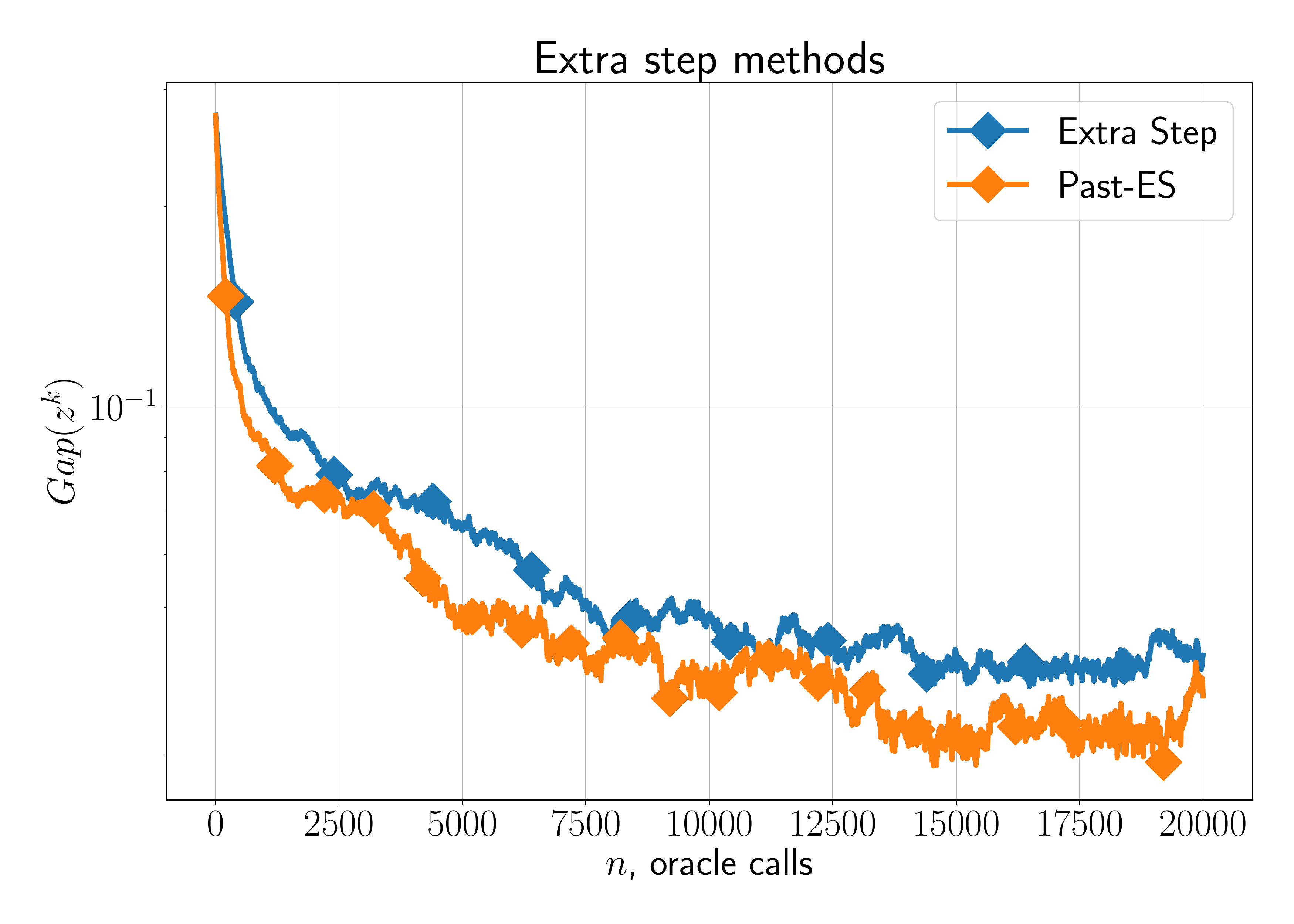}
         \caption{{\tt Extra Step} and {\tt Past-ES}}
     \end{subfigure}
    %  \hspace{0.1cm}
     \begin{subfigure}[H]{0.32\textwidth}
         \centering
         \includegraphics[width=\textwidth]{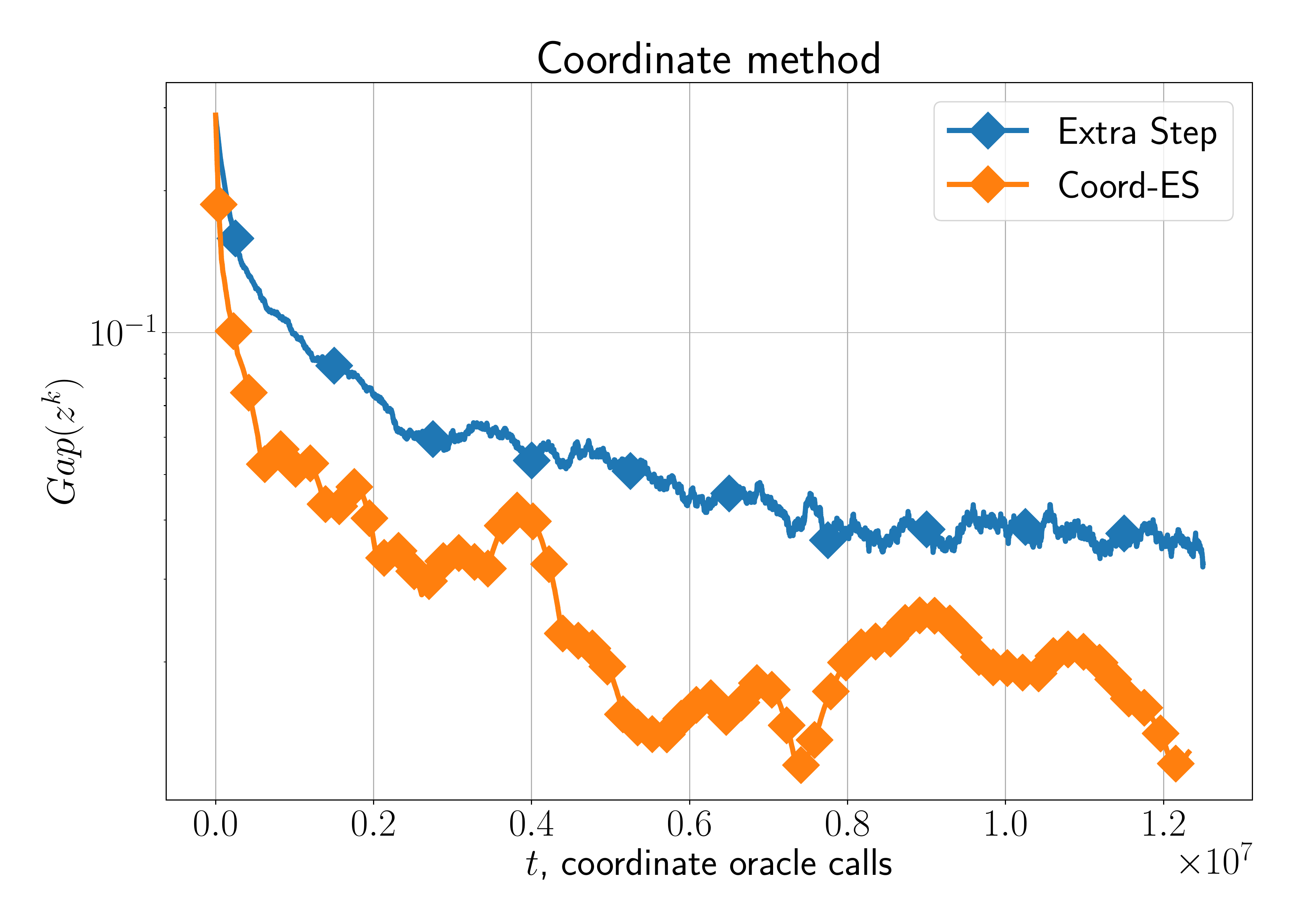}
         \caption{{\tt Extra Step} and {\tt Coord-ES}}
     \end{subfigure}
    %  \hspace{0.1cm}
     \begin{subfigure}[H]{0.32\textwidth}
         \centering
         \includegraphics[width=\textwidth]{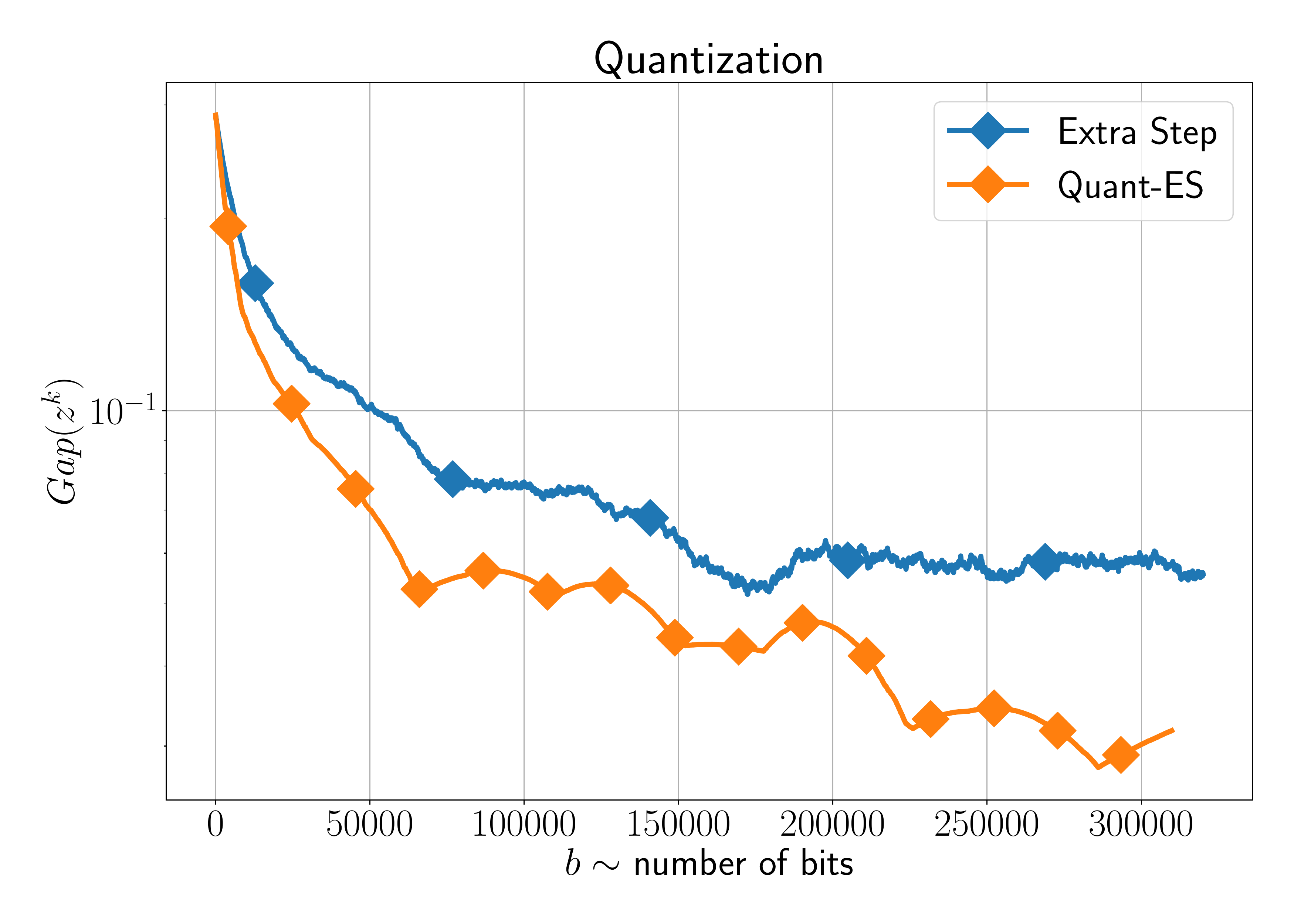}
         \caption{{\tt Extra Step} and {\tt Quant-ES}}
     \end{subfigure}
    \caption{Сходимость методов с дополнительным шагом для задачи  Полицейский vs. Грабитель \eqref{problem}.}
    \label{fig:numexp}
\end{figure}

\end{document}